\input ./style/arxiv-general.cfg
\documentclass[aos,MSNbibl,seceqn,dvips]{arximspdf}
\makeatletter
   \@ifpackageloaded{graphicx}{}{\usepackage{graphicx}}
\makeatother
\usepackage{accents}
\usepackage{url,breakurl}

\doi{10.1214/15-AOS1363}
\volume{44}
\issue{1}
\pubyear{2016}
\firstpage{183}
\lastpage{218}
\docsubty{FLA}

\makeatletter
\newcommand{\rrvert}{\vert}
\newcommand{\llvert}{\vert}
\renewcommand{\mathring}[1]{\accentset{\circ}{#1}}
\newtheorem{thmm}{Theorem}
\newtheorem{lma}{Lemma}
\newtheorem{cor}{Corollary}
\newtheorem{prop}{Proposition}

\makeatother

\begin{document}
\begin{frontmatter}

\title{Functional data analysis for density functions by transformation
to a Hilbert space}
\runtitle{Densities as functional data}

\begin{aug}
\author[A]{\fnms{Alexander} \snm{Petersen}\corref{}\ead[label=e1]{alxpetersen@gmail.com}}
\and
\author[A]{\fnms{Hans-Georg} \snm{M\"uller}\thanksref{T2}\ead[label=e2]{hgmueller@ucdavis.edu}}
\runauthor{A. Petersen and H.-G. M\"uller}
\affiliation{University of California, Davis}
\thankstext{T2}{Supported in part by NSF Grants DMS-11-04426,
DMS-12-28369 and DMS-14-07852.}
\address[A]{Department of Statistics\\
University of California, Davis\\
Mathematical Sciences Building 4118\\
399 Crocker Lane\\
One Shields Avenue\\
Davis, California 95616\\
USA\\
\printead{e1}\\
\phantom{E-mail:\ }\printead*{e2}}
\end{aug}

%
\received{\smonth{12} \syear{2014}}
%
\revised{\smonth{7} \syear{2015}}

%
\begin{abstract}
Functional data that are nonnegative and have a constrained integral
can be considered as samples of one-dimensional density functions. Such
data are ubiquitous. Due to the inherent constraints, densities do not
live in a vector space and, therefore, commonly used Hilbert space
based methods of functional data analysis are not applicable. To
address this problem, we introduce a transformation approach, mapping
probability densities to a Hilbert space of functions through a
continuous and invertible map. Basic methods of functional data
analysis, such as the construction of functional modes of variation,
functional regression or classification, are then implemented by using
representations of the densities in this linear space. Representations
of the densities themselves are obtained by applying the inverse map
from the linear functional space to the density space. Transformations
of interest include log quantile density and log hazard
transformations, among others. Rates of convergence are derived for the
representations that are obtained for a general class of
transformations under certain structural properties. If the
subject-specific densities need to be estimated from data, these rates
correspond to the optimal rates of convergence for density estimation.
The proposed methods are illustrated through simulations and
applications in brain imaging.
\end{abstract}

%
\begin{keyword}[class=AMS]
\kwd[Primary ]{62G05}
\kwd[; secondary ]{62G07}
\kwd{62G20}
\end{keyword}
\begin{keyword}
\kwd{Basis representation}
\kwd{kernel estimation}
\kwd{log hazard}
\kwd{prediction}
\kwd{quantiles}
\kwd{samples of density functions}
\kwd{rate of convergence}
\kwd{Wasserstein metric}
\end{keyword}
\end{frontmatter}

\section{Introduction}
\label{sec: Intro}

\begin{figure}

\includegraphics{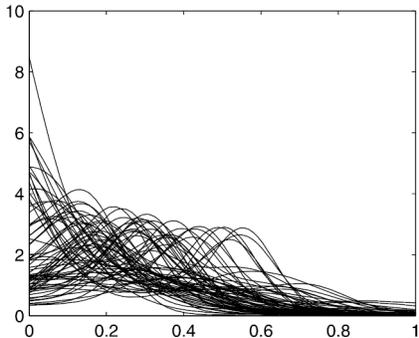}

\caption{Densities based on kernel density estimates for time course
correlations of BOLD signals obtained from brain fMRI between voxels in
a region of interest. Densities are shown for $n=68$ individuals
diagnosed with Alzheimer's disease. For details on density estimation,
see Section~\protect\ref{ss: dens_est}. Details regarding this data analysis,
which illustrates the proposed methods, can be found in
Section~\protect\ref{ss: brain}.}
\label{fig: rpl_cor}
\end{figure}

Data that consist of samples of one-dimensional distributions or
densities are common. Examples giving rise to such data are income
distributions for cities or states, distributions of the times when
bids are submitted in online auctions, distributions of movements in
longitudinal behavior tracking or distributions of voxel-to-voxel
correlations in fMRI signals (see Figure~\ref{fig: rpl_cor}). Densities
may also appear in functional regression models as predictors or responses.

The functional modeling of density functions is difficult due to the
two constrains \mbox{$\int f(x)\,dx=1$} and $f \ge0$. These
characteristics imply that the functional space where densities live is
convex but not linear, leading to problems for the application of
common techniques of functional data analysis (FDA) such as functional
principal components analysis (FPCA). This difficulty has been
recognized before and an approach based on compositional data methods
has been sketched in \cite{deli:11}, applying theoretical results in
\cite
{egoz:06}, which define a Hilbert structure on the space of densities.
Probably the first work on a functional approach for a sample of
densities is \cite{knei:01}, who utilized FPCA directly in density space
to analyze samples of time-varying densities and focused on the trends
of the functional principal components over time as well as the effects
of the preprocessing step of estimating the densities from actual
observations. Box--Cox transformations for a single nonrandom density
function were considered in \cite{wand:99}, who aimed at improving global
bandwidth choice for kernel estimation of a single density function.

Density functions also arise in the context of warping, or
registration, as time-warping functions correspond to distribution
functions. In the context of functional data and shape analysis, such
time-warping functions have been represented as square roots of the
corresponding densities \cite{sriv:07,sriv:11:02,sriv:11:01}, and these
square root densities reside in the Hilbert sphere, about which much is
known. For instance, one can define the Fr\'echet mean on the sphere
and also implement a nonlinear PCA method known as Principal Geodesic
Analysis (PGA) \cite{flet:04}. We will compare this alternative
methodology with our proposed approach in Section~\ref{sec: application}.

In this paper, we propose a novel and straightforward transformation
approach with the explicit goal of using established methods for
Hilbert space valued data once the densities have been transformed. The
key idea is to map probability densities into a linear function space
by using a suitably chosen continuous and invertible map $\psi$. Then
FDA methodology, which might range anywhere from exploratory techniques
to predictive modeling, can be implemented in this linear space. As an
example of the former, functional modes of variation can be constructed
by applying linear methods to the transformed densities, then mapping
back into the density space by means of the inverse map. Functional
regression or classification applications that involve densities as
predictors or responses are examples of the latter.

We also present theoretical results about the convergence of these
representations in density space under suitable structural properties
of the transformations. These results draw from known results for
estimation in FPCA and reflect the additional uncertainty introduced
through both the forward and inverse transformations. One rarely
observes data in the form of densities; rather, for each density, the
data are in the form of a random sample generated by the underlying
distribution. This fact will need to be taken into account for a
realistic theoretical analysis, adding a layer of complexity. Specific
examples of transformations that satisfy the requisite structural
assumptions are the log quantile density and the log hazard transformations.

A related approach can be found in a recent preprint by \cite{hron:14},
where the compositional approach of \cite{deli:11} was extended to define
a version of FPCA on samples of densities. The authors represent
densities by a centered log-ratio, which provides an isometric
isomorphism between the space of densities and the Hilbert space $L^2$,
and emphasize practical applications, but do not provide theoretical
support or consider the effects of density estimation. Our methodology
differs in that we consider a general class of transformations rather
than one specific transformation. In particular, the transformation can
be chosen independent of the metric used on the space of densities.
This provides flexibility since, for many commonly-used metrics on the
space of densities (see Section~\ref{ss: metrics}) corresponding
isometric isomorphisms do not exist with the $L^2$ distance in the
transformed space.

The paper is organized as follows: Pertinent results on density
estimation and background on metrics in density space can be found in
Section~\ref{sec: prelim}. Section~\ref{sec: fda} describes the basic
techniques of FPCA, along with their shortfalls when dealing with
density data. The main ideas for the proposed density transformation
approach are in Section~\ref{sec: transformation}, including an
analysis of specific transformations. Theory for this method is
discussed in Section~\ref{sec: theory}, with all proofs relegated to
the \hyperref[app]{Appendix}. In Section~\ref{ss: simulation}, we provide
simulations that illustrate the advantages of the transformation
approach over the direct functional analysis of density functions, also
including methods derived from properties of the Hilbert sphere. We
also demonstrate how densities can serve as predictors in a functional
regression analysis by using distributions of correlations of fMRI
brain imaging signals to predict cognitive performance. More details
about this application can be found in Section~\ref{ss: brain}.

\section{Preliminaries}
\label{sec: prelim}

\subsection{Density modeling}
\label{ss: model}

Assume that data consist of a sample of $n$ (random) density functions
$f_1, \ldots, f_n$, where the densities are supported on a common
interval $[0,T]$ for some $T>0$. Without loss of generality, we take $T
= 1$. The assumption of compact support is for convenience, and does
not usually present a problem in practice. Distributions with unbounded
support can be handled analogously if a suitable integration measure is
used. The main theoretical challenge for spaces of functions defined on
an unbounded interval is that the uniform norm is no longer weaker than
the $L^2$ norm, if the Lebesgue measure is used for the latter. This
can be easily addressed by replacing the Lebesgue measure $dx$ with a
weighted version, for example, $e^{-x^2} \,dx$.

Denote the space of continuous and strictly positive densities on $[0,
1]$ by $\mathcal {G}$. The sample consists of i.i.d. realizations of an
underlying stochastic process, that is, each density is independently
distributed as $f \sim\mathfrak{F}$, where $\mathfrak{F}$ is an $L^2$
process \cite{ash:75} on $[0,1]$ taking values in some space
$\mathcal {F}
\subset\mathcal {G}$. 
A basic assumption we make on the space $\mathcal {F}$ is:
\begin{longlist}[(A1)]
\item[(A1)] For all $f \in\mathcal {F}$, $f$ is continuously differentiable.
Moreover, there is a constant $M > 1$ such that, for all $f \in
\mathcal
{F}$, $\Vert f \Vert_\infty$, $\Vert1/f \Vert_\infty$ and $\Vert
f' \Vert_\infty$ are all bounded above
by $M$.
\end{longlist}

Densities $f$ can equivalently be represented as cumulative
distribution functions (c.d.f.) $F$ with domain $[0,1]$, hazard functions
$h=f/(1-F)$ (possibly on a subdomain of $[0,1]$ where $F(x)<1$) and
quantile functions $Q=F^{-1}$, with support $[0,1]$. Occasionally of
interest is the equivalent notion of the quantile-density function
$q(t)=Q'(t)=\frac{d}{dt}F^{-1}(t)= [f(Q(t)) ]^{-1}$, from which
we obtain $f(x)= [q(F(x)) ]^{-1}$, where we use the notation of
\cite{jone:92:1}. This concept goes back to \cite{parz:79} and \cite
{tuke:65}. Another classical notion of interest is the density-quantile
function $f(Q(t))$, which can be interpreted as a time-synchronized
version of the density function \cite{mull:11}. All of these functions
provide equivalent characterizations of distributions.

In many situations, the densities themselves will not be directly
observed. Instead, for each $i$, we may observe an i.i.d. sample of
data $W_{il}$, $l = 1, \ldots, N_i$, that are generated by the random
density $f_i$. Thus, there are two random mechanisms at work that are
assumed to be independent: the first generates the sample of densities
and the second generates the samples of real-valued random
data; one sample for each random density in the sample of densities.
Hence, the probability space can be thought of as a product space
$(\Omega_1\times\Omega_2, \mathcal{A}, P)$, where $P = P_1 \otimes P_2$.

\subsection{Metrics in the space of density functions}
\label{ss: metrics}

Many metrics and semimetrics on the space of density functions have
been considered, including the $L^2$, $L^1$ \cite{devr:85}, Hellinger and
Kullback--Leibler metrics, to name a few. In previous applied and
methodological work \cite{bols:03,mall:72,mull:11}, it was found that a
metric $d_{Q}$ based on quantile functions ${d_{Q}(f,g)^2=\int_0^1
(F^{-1}(t)-G^{-1}(t))^2 \,dt}$ is particularly promising from a
practical point of view.

This quantile metric has connections to the optimal transport problem
\cite{vill:03}, and corresponds to the Wasserstein metric between two
probability measures,
%
\begin{equation}
\label{eq: wass_dist} d_W(f,g)^2 = \inf_{X \sim f, Y \sim g}
E(X - Y)^2,
\end{equation}
where the expectation is with respect to the joint distribution of $(X,
Y)$. The equivalence $d_{Q}=d_{W}$ can be most easily seen by applying
a covariance identity due to~\cite{hoef:40}; details can be found in the
supplemental article \cite{pete:15}. We will develop our methodology for
a general metric, which will be denoted by $d$ in the following, and
may stand for any of the above metrics in the space of densities.
\subsection{Density estimation}
\label{ss: dens_est}

A common occurrence in functional data analysis is that the functional
data objects of interest are not completely observed. In the case of a
sample of densities, the information about a specific density in the
sample usually is available only through a random sample that is
generated by this density. Hence, the densities themselves must first
be estimated. Consider the estimation of a density $f \in\mathcal
{F}$ from
an i.i.d. sample (generated by $f$) of size $N$ by an estimator
$\check{f}$.
Here, $N = N(n)$ will implicitly represent a sequence that depends on
$n$, the size of the sample of random densities. In practice, any
reasonable estimator can be used that produces density estimates that
are bona fide densities and which can then be transformed into a linear
space. For the theoretical results reported in Section~\ref{sec:
theory}, a density estimator $\check{f}$ must satisfy the following
consistency properties in terms of the $L^2$ and uniform metrics
(denoted as $d_2$ and $d_{\infty}$, resp.):

\begin{longlist}[(D1)]
\item[(D1)] For a sequence $b_N = o(1)$, the density estimator $\check{f}$,
based on an i.i.d. sample of size $N$, satisfies $\check{f}\geq0$,
$\int_0^1
\check{f}
(x) \,dx = 1$ and
\[
\sup_{f\in\mathcal {F}} E\bigl(d_2(f, \check{f})^2
\bigr) = O\bigl(b_N^2\bigr).
\]
\item[(D2)] For a sequence $a_N = o(1)$ and some $R > 0$, the density
estimator $\check{f}$, based on an i.i.d. sample of size $N$, satisfies
\[
\sup_{f \in\mathcal {F}}P\bigl(d_{\infty}(f, \check{f}) >
Ra_N\bigr) \rightarrow 0.
\]
\end{longlist}

When this density estimation step is performed for densities on a
compact interval, which is the case in our current framework, the
standard kernel density estimator does not satisfy these assumptions,
due to boundary effects. Much work has been devoted to rectify the
boundary effects when estimating densities with compact support \cite
{cowl:96,mull:99}, but the resulting estimators leave the density space
and have not been shown to satisfy (D1) and (D2). Therefore, we
introduce here a
modified density estimator of kernel type that is guaranteed to satisfy
(D1) and (D2).

Let $\kappa $ be a kernel that corresponds to a continuous probability
density function and $h < 1/2$ be the bandwidth. We define a new kernel
density estimator to estimate the density $f \in\mathcal {F}$ on $[0,1]$
from a sample $W_1,\ldots,W_N \stackrel {\mathrm{i.i.d.}}{\sim} f$ by
%
\begin{equation}
\label{eq: dens_est} \check{f}(x) = \sum_{l = 1}^N
\kappa \biggl(\frac{x - W_l}{h} \biggr)w(x, h) \bigg/ \sum
_{l = 1}^N \int_0^1
\kappa \biggl(\frac{y - W_l}{h} \biggr)w(y, h) \,dy,
\end{equation}
for $x \in[0,1]$ and $0$ elsewhere. Here, the kernel $\kappa $ is assumed
to satisfy the following additional conditions:
\begin{longlist}[(K1)]
\item[(K1)] The kernel $\kappa $ is of bounded variation and is symmetric
about $0$.
\item[(K2)] The kernel $\kappa $ satisfies $\int_0^1 \kappa (u) \,du
> 0$,
and $\int_\mathbb{R} |u|\kappa (u) \,du$, $\int_\mathbb{R} \kappa
^2(u) \,du$
and $\int_\mathbb{R} |u|\kappa ^2(u) \,du$ are finite.
\end{longlist}
The weight function
\[
w(x, h) = \cases{ %
\displaystyle \biggl(\int_{-x/h}^1
\kappa (u) \,du \biggr)^{-1}, & \quad $\mbox{for } x\in[0,h),$
\vspace*{2pt}\cr
\displaystyle\biggl(\int_{-1}^{(1-x)/h} \kappa (u) \,du
\biggr)^{-1}, &\quad  $\mbox{for } x\in (1-h,1], \mbox{ and}$
\vspace*{2pt}\cr
1, & \quad $\mbox{otherwise},$ }
\]
is designed to remove boundary bias.

The following result 
demonstrates that this modified kernel estimator indeed satisfies
conditions (D1) and (D2). Furthermore, this result provides the rate in
(D1) for this estimator as $b_N = N^{-1/3}$, which is known to be the
optimal rate under our assumptions \cite{tsyb:09}, where the class of
densities $\mathcal {F}$ is assumed to be continuously differentiable,
and it
also shows that rates $a_N = N^{-c}$, for any $c \in(0, 1/6)$ are
possible in (D2).

\begin{prop}
\label{prop: dens_est}
If assumptions \textup{(A1)}, \textup{(K1)} and \textup{(K2)} hold, then the modified kernel
density estimator (\ref{eq: dens_est}) satisfies assumption \textup{(D1)}
whenever $h \rightarrow 0$ and $Nh \rightarrow \infty$ as $N
\rightarrow \infty$ with $b_N^2 =
h^2 + (Nh)^{-1}$. By taking $h = N^{-1/3}$ and $a_N = N^{-c}$ for any $c
\in(0, 1/6)$, \textup{(D2)} is also satisfied. In \textup{(S1)}, we may take $m(n) =
n^r$ for any $r > 0$.
\end{prop}

%
%

Alternative density estimators could also be used. In particular, the
beta kernel density estimator proposed in \cite{chen:99} is a promising
prospect. The convergence of the expected squared $L^2$ metric was
established in \cite{chen:99}, while weak uniform consistency was proved
in \cite{boue:03}. This density estimator is nonnegative, but requires
additional normalization to guarantee that it resides in the density space.

\section{Functional data analysis for the density process}
\label{sec: fda}

For a generic density function process $f \sim\mathfrak{F}$, denote
the mean function by $\mu(x) = E(f(x))$, the covariance function by
$G(x,y) = \operatorname{Cov}(f(x), f(y))$, and the orthonormal
eigenfunctions and
eigenvalues of the linear covariance operator $(Af)(t)=\int G(s,t)
f(s) \,ds$ by $\{\phi_k\}_{k = 1}^ \infty$ and $\{\lambda_k\}_{k = 1}^
\infty$, where the latter are positive and in decreasing order. If
$f_1, \ldots, f_n$ are i.i.d. distributed as $f$, then by the
Karhunen--Lo\`eve expansion, for each $i$,
%
\[
f_i(x) = \mu(x) + \sum_{k = 1}^ \infty
\xi_{ik} \phi_k(x),
\]
%
where $\xi_{ik} = \int_0^1 (f_i(x) - \mu(x))\phi_k(x) \,dx$ are the
uncorrelated principal components with zero mean and variance $\lambda
_k$. The Karhunen--Lo\`eve expansion constitutes the foundation for the
commonly used FPCA technique \cite
{bali:11,benk:06,bess:86,daux:82,hall:06:1,mull:06:7,li:10}.

The mean function $\mu$ of a density process $\mathfrak{F}$ is also a density
function, as the space of densities is convex, and can be estimated
by
%
\[
\tilde{\mu}(x) = \frac{1}{n}\sum_{i=1}^nf_i(x)\quad
\mbox{respectively}\quad \hat{\mu}(x) = \frac{1}{n}\sum
_{i=1}^n\check{f}_i(x),
\]
%
where the version $\tilde{\mu}$ corresponds to the case when the
densities are
fully observed and the version $\hat{\mu}$ corresponds to the case
when they
are estimated using suitable estimators such as (\ref{eq: dens_est});
this distinction will be used throughout. However, in the common
situation where one encounters horizontal variation in the densities,
this mean is not a good measure of center. This is because the
cross-sectional mean can only capture vertical variation. When
horizontal variation is present, the $L^2$ metric does not induce an
adequate geometry on the density space. A better method is \emph
{quantile synchronization} \cite{mull:11}, a version of which has been
introduced in \cite{bols:03} in the context of a genomics application.
Essentially, this involves considering the cross-sectional mean
function, $Q_\oplus (t)=E(Q(t))$, of the corresponding\vspace*{1pt} quantile process,
$Q$. The synchronized mean density is then given by $f_\oplus =
(Q_\oplus
^{-1})'$.

The quantile synchronized mean can be interpreted as a Fr\'echet mean
with respect to the Wasserstein metric $d=d_W$, where for a metric $d$
on $\mathcal {F}$ the Fr\'echet mean of the process $\mathfrak{F}$ is
defined by
%
\begin{equation}
\label{eq: frechet} f_\oplus = \arg\inf_{g \in\mathcal {F}} E\bigl(d(f,
g)^2\bigr),
\end{equation}
and the Fr\'echet variance is $E(d(f, f_\oplus )^2)$. Hence, for the
choice $d = d_W$, the Fr\'echet mean coincides with the quantile
synchronized mean. Further discussion of this Wasserstein--Fr\'echet
mean and its estimation is provided in the supplemental article \cite
{pete:15}. Noting that the cross-sectional mean corresponds to the Fr\'
echet mean for the choice $d = d_2$, the Fr\'echet mean provides a
natural measure of center, adapting to the chosen metric or geometry.

Modes of variation \cite{cast:86} have proved particularly useful in
applications to interpret and visualize the Karhunen--Lo\`eve
representation and FPCA \cite{jone:92,rams:05}. They focus on the
contribution of each eigenfunction $\phi_k$ to the stochastic behavior
of the process. The $k$th mode of variation is a set of functions
indexed by a parameter $\alpha\in\mathbb {R}$ that is given by
%
\begin{equation}
\label{eq: modes} g_k(x, \alpha) = \mu(x) + \alpha\sqrt{
\lambda_k}\phi_k(x).
\end{equation}
In order to construct estimates of these modes, and generally to
perform FPCA, the following estimates of the covariance function $G$ of
$\mathfrak{F}$ are needed:
\begin{eqnarray*}
\widetilde {G}(x,y) &= &\frac{1}{n}\sum_{i=1}^nf_i(x)f_i(y)
- \tilde {\mu}(x)\tilde{\mu}(y) \quad\mbox {respectively}
\\
\widehat {G}(x,y) &=& \frac{1}{n}\sum_{i=1}^n
\check{f}_i(x)\check {f}_i(y) - \hat{\mu}(x)\hat{
\mu}(y). %
\end{eqnarray*}
The eigenfunctions of the corresponding covariance operators, $\tilde
{\phi}_k$
or $\hat{\phi}_k$, then serve as estimates of $\phi_k$. Similarly, the
eigenvalues $\lambda_k$ are estimated by the empirical eigenvalues
($\tilde{\lambda}_k$ or $\hat{\lambda}_k$).

The empirical modes of variation are obtained by substituting estimates
for the unknown quantities in the modes of variation (\ref{eq: modes}),
%
\[
\tilde{g}_k(x, \alpha) = \tilde{\mu}(x) + \alpha\sqrt{\tilde {
\lambda}_k}\tilde{\phi}_k(x) \quad\mbox{respectively}\quad
\hat{g}_k(x, \alpha) = \hat{\mu}(x) + \alpha\sqrt{\hat{\lambda
}_k}\hat{\phi} _k(x).%
\]
%
These modes are useful for visualizing the FPCA in a Hilbert space. In
a nonlinear space such as the space of densities, they turn out to be
much less useful. Consider the eigenfunctions $\phi_k$. In \cite
{knei:01}, it was observed that estimates of these eigenfunctions for
samples of densities satisfy $\int_0^1 \hat{\phi}_k(x) \,dx = 0$ for
all $k$.
Indeed, this is true of the population eigenfunctions as well. To see
this, consider the following argument. Let ${\mathbf 1}(x) \equiv1$ so
that $\langle f - \mu, {\mathbf 1}\rangle= 0$. Take $\varphi$ to be the
projection of $\phi_1$ onto $\{{\mathbf 1}\}^\perp$. It is clear that
$\Vert\varphi\Vert_2 \leq1$ and $\operatorname{Var}(\langle f -
\mu, \phi_1
\rangle) = \operatorname{Var}(\langle f - \mu, \varphi\rangle)$.
However, by
definition, $ \operatorname{Var}(\langle f - \mu, \phi_1\rangle) =
\max_{\Vert
\phi
\Vert_2 = 1} \operatorname{Var}(\langle f - \mu, \phi\rangle)$.
Hence, in order to
avoid a contradiction, we must have $\Vert\varphi\Vert_2 = 1$, so
that $\langle\phi_1, {\mathbf 1} \rangle= 0$. The proof for all of the
eigenfunctions follows by induction.

At first, this seems like a desirable characteristic of the
eigenfunctions since it enforces \mbox{$\int g_k(x, \alpha) \,dx = 1$}
for any $k$ and $\alpha$. However, for $|\alpha|$ large enough, the
resulting modes of variation leave the density space since $\langle
\phi
_k, 1 \rangle=0$ implies at least one sign change for all
eigenfunctions. This also has the unfortunate consequence that the
modes of variation intersect at a fixed point which, as we will see in
Section~\ref{sec: application}, is an undesirable feature for
describing variation of samples of densities.

In practical applications, it is customary to adopt a
finite-dimensional approximation of the random functions by a truncated
Karhunen--Lo\`eve representation, including the first $K$ expansion terms,
%
\begin{equation}
\label{eq: fnt_trnc} f_i(x, K) = \mu(x) + \sum
_{k = 1}^K \xi_{ik}\phi_k(x).
\end{equation}
Then the functional principal components (FPC) $\xi_{ik}, k=1,\ldots
, K$, are used to represent each sample function. For fully observed
densities, estimates of the FPCs are obtained through their
interpretation as inner products,
%
\[
\tilde{\xi}_{ik} = \int_0^1
\bigl(f_i(x) - \tilde{\mu}(x)\bigr)\tilde{\phi }_k(x)
\,dx.
\]
%
The truncated processes in (\ref{eq: fnt_trnc}) are then estimated by
simple plug-in. Since the truncated finite-dimensional representations
as derived from the finite-dimensional Karhunen--Lo\`eve expansion are
designed for functions in a linear space, they are good approximations
in the $L^2$ sense, but (i) may lack the defining characteristics of a
density and (ii) may not be good approximations in a nonlinear space.

Thus, while it is possible to directly apply FPCA to a sample of
densities, this approach provides an extrinsic analysis as the ensuing
modes of variation and finite-dimensional representations leave the
density space. One possible remedy would be to project these quantities
back onto the space of densities, say by taking the positive part and
renormalizing. In the applications presented in Section~\ref{sec:
application}, we compare this ad hoc procedure with the proposed
transformation approach.

\section{Transformation approach}
\label{sec: transformation}

The proposed transformation approach is to map the densities into a new
space $L^2(\mathcal {T})$ via a functional transformation $\psi$, where
$\mathcal
{T} \subset\mathbb {R}$ is a compact interval. Then we work with the
resulting $L^2$ process $X: = \psi(f)$. By performing FPCA in the
linear space $L^2(\mathcal {T})$ and then mapping back to density
space, this
transformation approach can be viewed as an intrinsic analysis, as
opposed to ordinary FPCA. With $\nu$ and $H$ denoting the mean and
covariance functions, respectively, of the process $X$, $\{\rho_k\}_{k
= 1}^\infty$ denoting the orthonormal eigenfunctions of the covariance
operator with kernel $H$ with corresponding eigenvalues $\{\tau_k\}_{k
= 1}^\infty$, the Karhunen--Lo\`eve expansion for each of the
transformed processes $X_i = \psi(f_i)$ is
%
\[
X_i(t) = \nu(t) + \sum_{k = 1}^\infty
\eta_{ik}\rho_k(t), \qquad t \in \mathcal {T},
\]
%
with principal components $\eta_{ik} = \int_\mathcal {T} (X_i(t) -
\nu
(t))\rho
_k(t) \,dt$.

Our goal is to find suitable transformations $\psi: \mathcal {G}
\rightarrow
L^2(\mathcal {T})$ from density space to a linear functional space. To be
useful in practice and to enable derivation of consistency properties,
the maps $\psi$ and $\psi^{-1}$ must satisfy certain continuity
requirements, which will be given at the end of this section. We begin
with two specific examples of relevant transformations. For clarity,
for functions in the native density space $\mathcal {G}$ we denote the
argument by $x$, while for functions in the transformed space
$L^2(\mathcal
{T})$ the argument is $t$.

\textit{The log hazard transformation.} Since hazard functions diverge at
the right endpoint of the distribution, which is 1, we consider
quotient spaces induced by identifying densities which are equal on a
subdomain $\mathcal {T} = [0, {1_\delta}]$, where ${1_\delta}= 1 -
\delta$ for
some $0 <
\delta< 1$. With a slight abuse of notation, we denote this quotient
space as $\mathcal {G}$ as well. The log hazard transformation $\psi
_H: \mathcal
{G} \rightarrow L^2(\mathcal {T})$ is
\[
\psi_H(f) (t)= \log\bigl(h(t)\bigr) = \log \biggl\{
\frac{f(t)}{1-F(t)} \biggr\}, \qquad t \in\mathcal {T}.
\]
Since the hazard function is positive but otherwise not constrained on
$\mathcal {T}$, it is easy to see that $\psi$ indeed maps density functions
to $L^2(\mathcal {T})$.
The inverse map can be defined for any continuous function $X$ as
\[
\psi_H^{-1}(X) (x) = \exp \biggl\{X(x) - \int
_0^x e^{X(s)} \,ds \biggr\} ,\qquad  x \in[0,
{1_\delta}].
\]
Note that for this case one has a strict inverse only modulo the
quotient space. However, in order to use metrics such as $d_W$, we must
choose a representative. A~straightforward way to do this is to assign
the remaining mass uniformly, that is,
\[
\psi_H^{-1}(X) (x) = \delta^{-1}\exp \biggl\{ -
\int_0^{{1_\delta}} e^{X(s)} \,ds \biggr\},\qquad  x
\in({1_\delta}, 1].
\]

\textit{The log quantile density transformation.} For $\mathcal {T} = [0,1]
$, the
log quantile density (LQD) transformation $\psi_Q: \mathcal {G}
\rightarrow L^2(\mathcal
{T})$ is given by
\[
\psi_Q(f) (t) = \log\bigl(q(t)\bigr) = -\log\bigl\{f\bigl(Q(t)
\bigr)\bigr\}, \qquad t \in\mathcal {T}.
\]
It is then natural to define the inverse of a continuous function $X$
on $\mathcal {T}$ as the density given by
$\exp\{-X(F(x))\}$, where $Q(t)=F^{-1}(t) = \int_0^t e^{X(s)} \,ds$.
Since the value $F^{-1}(1)$ is not fixed, the support of the densities
is not fixed within the transformed space,
and as the inverse transformation should map back into the space of
densities with support on $[0, 1]$, we make a slight adjustment when
defining the inverse by
\[
\psi_Q^{-1}(X) (x) = \theta_X\exp\bigl\{-X
\bigl(F(x)\bigr)\bigr\},\qquad  F^{-1}(t) = \theta _X^{-1}
\int_0^t e^{X(s)} \,ds,
\]
where $\theta_X = \int_0^1 e^{X(s)} \,ds$. Since $F^{-1}(1) = 1$
whenever $X
\in\psi_Q (\mathcal {G} )$, this definition coincides
with the
natural definition mentioned above on $\psi_Q (\mathcal
{G} )$.

To avoid the problems that afflict the linear-based modes of variation
as described in Section~\ref{sec: fda}, in the transformation approach we construct
modes of variation in the transformed space for processes $X = \psi(f)$
and then map these back into the density space, defining transformation
modes of variation
%
\begin{equation}
\label{eq: tmodes} g_k(x, \alpha, \psi) = \psi^{-1} (\nu+
\alpha\sqrt{\tau _k}\rho _k ) (x).
\end{equation}
Estimation of these modes is done by first estimating the mean function
$\nu$ and covariance function $H$ of the process $X$. Letting
$\widehat {X}_i =
\psi(\check{f}_i)$, the empirical estimators are
%
\begin{eqnarray}
\tilde{\nu}(t) &= &\frac{1}{n}\sum_{i=1}^nX_i(t)\quad
\mbox{respectively}\quad  \hat{\nu}(t) = \frac{1}{n}\sum
_{i=1}^n\widehat {X}_i(t);
\label
{eq: nu_est}
\\
\label{eq: H_est}
\widetilde {H}(s,t) &=& \frac{1}{n}\sum_{i=1}^nX_i(s)X_i(t)
- \tilde {\nu}(s)\tilde{\nu}(t)\quad \mbox {respectively}
\nonumber
\\[-8pt]
\\[-8pt]
\nonumber
\widehat {H}(s,t) &=& \frac{1}{n}\sum_{i=1}^n
\widehat {X}_i(s)\widehat {X}_i(t) - \hat{\nu}(s)\hat{
\nu} (t).
\end{eqnarray}
Estimated eigenvalues and eigenfunctions ($\tilde{\tau}_k$ and
$\tilde{\rho}_k$,
resp., $\hat{\tau}_k$ and $\hat{\rho}_k$) are then obtained
from the mean
and covariance estimates as before, yielding the transformation mode of
variation estimators
%
\begin{eqnarray}\label{eq: tmodes_est}
\tilde{g}_k(x, \alpha, \psi) &=& \psi^{-1}(\tilde{\nu}+
\alpha \sqrt {\tilde{\tau}_k}\tilde{\rho} _k) (x)\quad
\mbox{respectively}
\nonumber
\\[-8pt]
\\[-8pt]
\nonumber
\hat{g}_k(x, \alpha, \psi) &= &\psi^{-1}(\hat{\nu}+ \alpha
\sqrt {\hat{\tau} _k}\hat{\rho} _k) (x).
\nonumber
\end{eqnarray}

In contrast to the modes of variation resulting from ordinary FPCA in
(\ref{eq: modes}), the transformation modes are bona fide density
functions for any value of $\alpha$. Thus, for reasonably chosen
transformations, the transformation modes can be expected to provide a
more interpretable description of the variability contained in the
sample of densities. Indeed, the data application in Section~\ref{ss:
brain} shows that this is the case, using the log quantile density
transformation as an example.

The truncated representations of the original densities in the sample
are then given by
%
\begin{equation}
\label{eq: fnt_trnc1} f_i(x, K, \psi) = \psi^{-1} \Biggl(\nu+
\sum_{k = 1}^K \eta _{ik}\rho
_k \Biggr) (x).
\end{equation}
Utilizing (\ref{eq: nu_est}), (\ref{eq: H_est}) and the ensuing
estimates of the eigenfunctions, the (transformation) principal
components, for the case of fully observed densities, are obtained in a
straightforward manner,
%
\begin{equation}
\label{eq: pc_est} \tilde{\eta}_{ik}= \int_\mathcal {T}
\bigl(X_i(t) - \tilde{\nu }(t)\bigr)\tilde{\rho}_k(t)
\,dt,
\end{equation}
whence
%
\[
\tilde{f}_i(x, K, \psi) = \psi^{-1} \Biggl(\tilde{\nu}+
\sum_{k =
1}^K \tilde{\eta}_{ik}
\tilde{\rho} _k \Biggr) (x).
\]
%

In practice, the truncation point $K$ can be selected by choosing a
cutoff for the fraction of variance explained. This raises the question
of how to quantify \emph{total variance}. For the chosen metric $d$, we
propose to use the Fr\'echet variance
%
\begin{equation}
\label{eq: tot_var} V_\infty:= E\bigl(d(f, f_\oplus )^2
\bigr),
\end{equation}
which is estimated by its empirical version
%
\begin{equation}
\label{eq: tot_var_est} \tilde{V}_\infty = \frac{1}{n}\sum
_{i=1}^nd(f_i, \tilde{f}_\oplus
)^2,
\end{equation}
using an estimator $\tilde{f}_\oplus $ of the Fr\'echet mean.
Truncating at $K$
included components as in (\ref{eq: fnt_trnc}) or in (\ref{eq:
fnt_trnc1}) and denoting the truncated versions as $f_{i,K}$,
the variance explained by the first $K$ components is
%
\begin{equation}
\label{eq: var_ex} V_K:= V_\infty - E\bigl(d(f_1,
f_{1, K})^2\bigr),
\end{equation}
which is estimated by
%
\begin{equation}
\label{eq: var_ex_est} \tilde{V}_K= \tilde{V}_\infty -
\frac{1}{n}\sum_{i=1}^nd(f_i,
\tilde{f}_{i, K})^2.
\end{equation}
The ratio $V_K/V_\infty $ is called the fraction of variance explained
(FVE), and is estimated by $\tilde{V}_K/\tilde{V}_\infty $. If the
truncation level is
chosen so that a fraction $p$, $0 < p < 1$, of total variation is to be
explained,
the optimal choice of $K$ is
%
\begin{equation}
\label{eq: choice_K} K^\ast= \min \biggl\{K : \frac{V_K}{V_\infty} > p \biggr
\},
\end{equation}
which is estimated by
%
\begin{equation}
\label{eq: choice_K_est} \tilde{K}^\ast= \min \biggl\{K : \frac{\tilde{V}_K}{\tilde
{V}_\infty} > p
\biggr\}.
\end{equation}
As will be demonstrated in the data illustrations, this more general
notion of variance explained is a useful concept when dealing with
densities or other functions that are not in a Hilbert space.
Specifically, we will show that density representations in (\ref{eq:
fnt_trnc1}), obtained via transformation, yield higher FVE values than
the ordinary representations in (\ref{eq: fnt_trnc}), thus giving more
efficient representations of the sample of densities.

For the theoretical analysis of the transformation approach, certain
structural assumptions on the transformations
need to be satisfied. The required smoothness properties for maps $\psi
$ and $\psi^{-1}$ are implied by the three conditions (T0)--(T3) below.
Here, the $L^2$ and uniform metrics are denoted by $d_2$ and $d_{\infty}$,
respectively, and the uniform norm is denoted by $\Vert\cdot \Vert
_\infty$.
\begin{longlist}[(T0)]
\item[(T0)] Let $f$, $g \in\mathcal {G}$ with $f$ differentiable and
$\Vert f' \Vert_\infty < \infty$. Set
\[
D_0 \geq\max \bigl(\Vert f \Vert_\infty, \Vert1/f
\Vert_\infty, \Vert g \Vert_\infty, \Vert1/g \Vert_\infty,
\bigl\Vert f' \bigr\Vert_\infty \bigr).
\]
Then there exists $C_0$ depending only on $D_0$ such that
\[
d_2\bigl(\psi(f), \psi(g)\bigr) \leq C_0
d_2(f, g),\qquad  d_{\infty}\bigl(\psi(f), \psi (g)\bigr) \leq
C_0 d_{\infty}(f, g).
\]
\item[(T1)] Let $f \in\mathcal {G}$ be differentiable with $\Vert f'
\Vert_\infty <
\infty$ and let $D_1$ be a constant bounded below by $\max
(\Vert f \Vert_\infty, \Vert1/f \Vert_\infty, \Vert f' \Vert
_\infty )$. Then $\psi(f)$ is differentiable
and there exists $C_1 > 0 $ depending only on $D_1$ such that $\Vert
\psi(f) \Vert_\infty \leq C_1$ and $\Vert\psi(f)' \Vert_\infty
\leq C_1$.
\item[(T2)] Let $d$ be the selected metric in density space, $Y$ be
continuous and $X$ be differentiable on $\mathcal {T}$ with $\Vert X'
\Vert_\infty<\infty
$. There exist constants $C_2 = C_2(\Vert X \Vert_\infty,
\Vert X' \Vert_\infty) > 0$
and $C_3 = C_3(d_{\infty}(X, Y)) > 0$ such that
\[
d\bigl(\psi^{-1}(X), \psi^{-1}(Y)\bigr) \leq C_2
C_3 d_2(X, Y)
\]
and, as functions, $C_2$ and $C_3$ are increasing in their respective arguments.
\item[(T3)] For a given metric $d$ on the space of densities and $f_{1,
K} = f_1(\cdot, K, \psi)$ [see~(\ref{eq: fnt_trnc1})], $V_\infty - V_K
\rightarrow
0$ and \mbox{$E(d(f, f_{1, K})^4) = O(1)$} as $K \rightarrow \infty$.
\end{longlist}

Here, assumptions (T0) and (T2) relate to the continuity of $\psi$ and
$\psi^{-1}$, while (T1) means that bounds on densities in the space
$\mathcal{G}$ are accompanied by corresponding bounds of the
transformed processes $X$. Assumption (T3) is needed to ensure that the
finitely truncated versions in the transformed space are consistent, as
the truncation parameter increases.

To establish these properties for the log hazard and log quantile
density transformations, 
denoting as before the mean function, covariance function,
eigenfunctions and eigenvalues associated with the process $X$ by $(\nu
, H, \rho_k, \tau_k)$, assumption (T1) implies that $\nu$, $H$,
$\rho
_k$, $\nu'$ and $\rho_k'$ are bounded for all $k$ (see Lemma~\ref{lma:
fpca_bounds} in the \hyperref[app]{Appendix} for details). In turn, these bounds
imply a nonrandom Lipschitz constant for the residual process $X - X_K
= \sum_{k = K + 1}^\infty\eta_{k}\phi_k$ as follows. Under (A1), the
constant $C_1$ in (T1) can be chosen uniformly over $f\in\mathcal
{F}$. As a
consequence, we have $\Vert X \Vert_\infty < C_1$ almost surely so
that $\Vert\nu \Vert_\infty
< C_1$ and
%
\begin{equation}\quad
\label{eq: fpc_bound} |\eta_k| = \biggl\llvert \int_{\mathcal {T}}
\bigl(X(t) - \nu(t)\bigr)\phi_k(t) \,dt\biggr\rrvert
\leq2C_1\int_{\mathcal {T}} \bigl|\phi_k(t)\bigr| \,dt
\leq2C_1|\mathcal {T}|^{1/2},
\end{equation}
almost surely. Additionally, $\Vert\nu' \Vert_\infty < C_1$ and
$\Vert\rho _k' \Vert_\infty <
\infty$ for all $k$ by dominated convergence, so that
\[
\bigl\Vert X_K' \bigr\Vert_\infty \leq\bigl\Vert
\nu'\bigr \Vert_\infty + \sum_{k =
1}^K
|\eta_k|\bigl\Vert\rho_k' \bigr\Vert_\infty
\leq C_1 \Biggl(1 + 2|\mathcal {T}|^{1/2}\sum
_{k = 1}^K\bigl\Vert\rho _k'
\bigr\Vert_\infty \Biggr).
\]
Since $\Vert X' \Vert_\infty < C_1$ almost surely, setting
%
\begin{equation}
\label{eq: lip_const} L_K:= 2C_1 \Biggl(1 + |\mathcal
{T}|^{1/2}\sum_{k = 1}^K \bigl\Vert\rho
_k' \bigr\Vert_\infty \Biggr)
\end{equation}
then yields the almost sure bound
\[
\bigl|(X - X_K) (s) - (X - X_K) (t)\bigr| \leq L_K|s -
t|.
\]

The following result demonstrates the continuity of the log hazard and
log quantile density transformations for classes of processes $X$ that
have suitably fast declining eigenvalues and suitable smoothness of the
finite approximations.

\begin{prop}
\label{prop: transformations}
Assumptions \textup{(T0)--(T2)} are satisfied for both $\psi_H$ and $\psi_Q$
with either $d = d_2$ or $d = d_W$. Let $L_K$ denote the Lipschitz
constant given in (\ref{eq: lip_const}). If:
\begin{longlist}[(ii)]
\item[(i)]$L_K\sum_{k = K+1}^\infty\tau_k = O(1)$ as $K \rightarrow
\infty$ and
\item[(ii)] there is a sequence $r_m$, $m \in\mathbb{N}$, such that $E(\eta
_{1k}^{2m}) \leq r_m\tau_k^m$ for large $k$ and $ (\frac
{r_{m+1}}{r_m} )^{1/3} = o(m)$,
\end{longlist}
are satisfied, then assumption \textup{(T3)} is also satisfied for both $\psi_H$
and $\psi_Q$ with either $d = d_2$ or $d = d_W$.
\end{prop}

As example, consider the Gaussian case for transformed processes $X$
[or, similarly, the truncated Gaussian case in light of (\ref{eq:
fpc_bound})] with components $\eta_{1k} \sim N(0,\lambda_k)$. Then
$E(\eta_{1k}^{2m})=\tau_k^m (2m-1)!!$, whence \mbox{$r_m=(2m-1)!!$} so
that \mbox{$ (r_{m+1}/{r_m} )^{1/3} = o(m)$} in (ii) is
trivially satisfied. If the eigenfunctions correspond to the
trigonometric basis, then $\Vert\rho_k' \Vert_\infty = O(k)$, so
that $L_K =
O(K^2)$. Hence, any eigenvalue sequence satisfying $\tau_k = O(k^{-4})$
would satisfy (i) in this case.

\section{Theoretical results}
\label{sec: theory}

The transformation modes of variation as defined in (\ref{eq: tmodes}),
together with the FVE values and optimal truncation points in~(\ref{eq:
choice_K}), constitute the main components of the proposed approach. In
this section, we investigate the weak consistency of the estimators of
these quantities, given in (\ref{eq: tmodes_est}) and (\ref{eq:
choice_K_est}), respectively, for the case of a generic density metric
$d$, as $n \rightarrow \infty$. While asymptotic properties of
estimates in
FPCA are well established \cite{bosq:00,li:10}, the effects of density
estimation and transformation need to be studied in order to validate
the proposed transformation approach. When densities are estimated, a
lower bound $m$ on the sample sizes available for estimating each
density is required, as stipulated in the following assumption:
\begin{longlist}[(S1)]
\item[(S1)] Let $\check{f}$ be a density estimator that satisfies
(D2), and
suppose densities $f_i \in\mathcal {F}$ are estimated by $\check
{f}_i$ from i.i.d.
samples of size $N_i = N_i(n)$, $i = 1,\ldots, n$, respectively. There
exists a sequence of lower bounds $m(n) \leq\break \min_{1 \leq i
\leq
n} N_i$ such that $m(n) \rightarrow \infty$ as $n \rightarrow
\infty$ and
\[
n\sup_{f \in\mathcal {F}}P\bigl(d_{\infty}(f, \check{f}) >
Ra_m\bigr) \rightarrow 0,
\]
where, for generic $f \in\mathcal {F}$, $\check{f}$ is the estimated
density from
a sample of size $N(n) \geq m(n)$.
\end{longlist}

Proposition~\ref{prop: dens_est} in Section~\ref{ss: dens_est} implies
that, for the density estimator in (\ref{eq: dens_est}), property (S1)
is satisfied for sequences of the form $m(n) = n^r$ for arbitrary $r >
0$. For $r < 3/2$, this rate dominates the rate of convergence in
Theorem~\ref{thmm: modes} below, which thus cannot be improved under our assumptions.
While the theory we provide is general in terms of the transformation
and metric, of particular interest are the specific transformations
discussed in Section~\ref{sec: transformation} and the Wasserstein
metric $d_W$. Proofs and auxiliary lemmas are in the \hyperref[app]{Appendix}.

To study the transformation modes of variation, auxiliary results
involving convergence of the mean, covariance, eigenvalue and
eigenfunction estimates in the transformed space are needed. These
auxiliary results are given in Lemma~\ref{lma: mean_cov} and Corollary~\ref{cor: fpca_est} in the \hyperref[app]{Appendix}. A critical component in
these rates is the spacing between eigenvalues
%
\begin{equation}
\delta_k = \min_{1 \leq j \leq k}(\tau_j -
\tau_{j + 1}). \label{eq:
sp}
\end{equation}
These spacings become important as one aims to estimate an increasing
number of transformation modes of variation simultaneously.

The following result provides the convergence of estimated
transformation modes of variation in (\ref{eq: tmodes_est}) to the true
modes $g_k(\cdot, \alpha, \psi)$ in (\ref{eq: tmodes}), uniformly over
mode parameters $|\alpha| \le\alpha_0$ for any constant $\alpha_0 >0$.
For the case of estimated densities, if (D1), (D2) and (S1) are
satisfied, $m=m(n)$ denotes the increasing sequence of lower bounds in
(S1), and $b_m$ is the rate of convergence in (D1), indexed by the
bounding sequence $m$.

\begin{thmm}
\label{thmm: modes}
Fix $K$ and $\alpha_0 > 0$. Under assumptions \textup{(A1)}, \textup{(T1)} and \textup{(T2)}, and
with $\tilde{g}_k, \hat{g}_k$ as in (\ref{eq: tmodes_est}),
\[
\max_{1 \leq k \leq K}\sup_{|\alpha| \leq\alpha_0}d\bigl(g_k(
\cdot, \alpha, \psi), \tilde{g}_k(\cdot, \alpha, \psi)\bigr) =
O_p\bigl(n^{-1/2}\bigr).
\]
Additionally, there exists a sequence $K(n) \rightarrow \infty$ such that
\[
\max_{1 \leq k \leq K(n)}\sup_{|\alpha| \leq\alpha_0}d\bigl(g_k(
\cdot, \alpha, \psi), \tilde{g}_k(\cdot, \alpha, \psi)\bigr) =
o_p(1).
\]
If assumptions \textup{(T0)}, \textup{(D1)}, \textup{(D2)} and \textup{(S1)} are also satisfied and $K$,
$\alpha_0$ are fixed,
\[
\max_{1 \leq k \leq K}\sup_{|\alpha| \leq\alpha_0}d\bigl(g_k(
\cdot, \alpha, \psi), \hat{g}_k(\cdot, \alpha, \psi)\bigr) =
O_p\bigl(n^{-1/2}+ b_m\bigr).
\]
Moreover, there exists a sequence $K(n) \rightarrow \infty$ such that
\[
\max_{1 \leq k \leq K(n)}\sup_{|\alpha| \leq\alpha_0}d\bigl(g_k(
\cdot, \alpha, \psi), \hat{g}_k(\cdot, \alpha, \psi)\bigr) =
o_p(1).
\]
\end{thmm}

In addition to demonstrating the convergence of the estimated
transformation modes of variation for both fully observed and estimated
densities, this result also provides uniform convergence over
increasing sequences of included components $K=K(n)$. Under assumptions
on the rate of decay of the eigenvalues and the upper bounds for the
eigenfunctions, one also can get rates for the case $K(n) \rightarrow
\infty$. For example, suppose the densities are fully observed, $\tau_k
= ce^{-\theta k}$ for \mbox{$c, \theta>0$} and $\sup_k\Vert\rho_k
\Vert_\infty
\leq A$ (as would be the case for the trigonometric basis, but this
could be easily replaced by a sequence $A_k$ of increasing bounds).
Additionally, suppose $C_2 = a_0e^{a_1\Vert X \Vert_\infty}$ in (T2),
as is the
case for the log quantile density transformation with the metric $d_W$
(see the proof of Proposition~\ref{prop: transformations}).\vspace*{1pt} Then,
following the proof of Theorem~\ref{thmm: modes}, one finds that, for
$K(n) = \lfloor\frac{1}{4\theta}\log n \rfloor$,
\[
{\max_{1 \leq k \leq K(n)}\sup_{|\alpha| \leq\alpha_0}d
\bigl(g_k(\cdot, \alpha, \psi), \tilde{g}_k(\cdot,
\alpha, \psi)\bigr) = O_p\bigl(n^{-1/4}\bigr)}.
\]

For the truncated representations in (\ref{eq: fnt_trnc1}), the
truncation point $K$ may be viewed as a tuning parameter. When adopting
the fraction of variance explained criterion [see (\ref{eq: tot_var})
and (\ref{eq: var_ex})]\vadjust{\goodbreak} for the data-adaptive selection of $K$, a user
will typically choose the fraction $p \in(0,1)$, for which the
corresponding optimal value $K^\ast$ is given in (\ref{eq: choice_K}),
with the data-based estimate in (\ref{eq: choice_K_est}). This requires
estimation of the Fr\'echet mean $f_\oplus $ (\ref{eq: frechet}), for
which we assume the availability of an estimator $\tilde{f}_\oplus $ that
satisfies $d(f_\oplus, \tilde{f}_\oplus ) = O_p(\gamma_n)$ for the
given metric $d$
in density space and some sequence $\gamma_n \rightarrow 0$. For the
choice $d
= d_W$, $\gamma_n = n^{-1/2}$ is admissible~\cite{pete:15}.

This selection procedure for the truncation parameter is a
generalization of the scree plot in multivariate analysis, where the
usual fraction of variance concept that is based on the eigenvalue
sequence is replaced here with the corresponding Fr\'echet variance. As
more data become available, it is usually desirable to increase the
fraction of variance explained in order to more accurately represent
the true underlying functions. Therefore, it makes sense to choose a
sequence $p_n \in(0,1)$, with $p_n \uparrow1$. The following result
provides consistent recovery of the fraction of variance explained
values $V_K/V_\infty $ as well as the optimal choice $K^\ast$ for
such sequences.

\begin{thmm}
\label{thmm: choice_K}
Assume \textup{(A1)} and \textup{(T1)--(T3)} hold. Additionally, suppose an estimator
$\tilde{f}_\oplus $ of $f_\oplus $ satisfies $d(f_\oplus , \tilde
{f}_\oplus ) = O_p(\gamma_n)$
for a sequence $\gamma_n \rightarrow0$. Then there is a sequence $p_n
\uparrow1$ such that
\[
\max_{1 \leq K \leq K^\ast}\biggl\llvert \frac{V_K}{V_\infty } -
\frac
{\tilde{V}_K
}{\tilde{V}_\infty
}\biggr\rrvert = o_p(1)
\]
and, consequently,
\[
P \bigl(K^\ast\neq\tilde{K}^\ast \bigr) \rightarrow 0.
\]
\end{thmm}

Specific choices for the sequence $p_n$ and their implications for the
corresponding sequence $K^\ast(n)$ can be investigated under additional
assumptions. For example, consider the case where $\tau_k =
ce^{-\theta
k}$, $\sup_k \Vert\rho_k \Vert_\infty\leq A$, ${V_\infty - V_K =
be^{-\omega K}}$,
$C_2 = a_0e^{a_1\Vert X \Vert_\infty}$ in (T2) and $\gamma_n =
n^{-1/2}$. Then, by
following the proofs of Lemma~\ref{lma: fnt_trnc} and Theorem~\ref{thmm:
choice_K}, we find that if $r<[2(2a_1C_1\times\break |\mathcal {T}|^{1/2}A + \theta+
\omega)]^{-1}$, the choice
\[
p_n = 1 - \frac{b(1 + e^\omega)}{2V_\infty }n^{-\omega r}
\]
leads to a corresponding sequence of tuning parameters $K^\ast(n) =
\lfloor r \log n\rfloor$. In particular, this means that
\[
\max_{1 \leq K \leq K^\ast}\biggl\llvert \frac{V_K}{V_\infty } -
\frac
{\tilde{V}_K
}{\tilde{V}_\infty
}\biggr\rrvert = O_p \biggl( \biggl(
\frac{\log n}{n} \biggr)^{1/2} \biggr)
\]
and the relative error $(\tilde{K}^\ast- K^\ast)/K^\ast$ converges
at the
rate $o_p(1/\log n)$ under these assumptions.

\section{Illustrations}
\label{sec: application}
\subsection{Simulation studies}
\label{ss: simulation}

Simulation studies were conducted to compare the performance between
ordinary FPCA applied to densities, the proposed transformation
approach using the log quantile density transformation, $\psi_Q$, and
methods derived for the Hilbert sphere \cite
{flet:04,sriv:07,sriv:11:02,sriv:11:01} for three simulation settings
that are listed in Table~\ref{tab: sim_set}. The first two settings
represent vertical and
horizontal variation, respectively, while the third setting is a
combination of both. We considered the case where the densities are
fully observed, as well as the more realistic case where only a random
sample of data generated by a density is available for each density. In
the latter case, densities were estimated from a sample of size 100
each, using the density estimator in (\ref{eq: dens_est}) with the
kernel $\kappa $ being the standard normal density and a bandwidth of $h
= 0.2$.
%
\begin{table}
\caption{Simulation designs for comparison of methods}\label{tab: sim_set}
\begin{tabular*}{\textwidth}{@{\extracolsep{\fill}}lcc@{}}
\hline
\textbf{Setting} & \textbf{Random component} & \textbf{Resulting density} \\
\hline
1 & $\log(\sigma_i)\sim\mathcal{U}[-1.5,1.5]$, $i = 1, \ldots,50$ &
$\mathcal{N}(0, \sigma_i^2)$ truncated on $[-3,3]$ \\
2 & $\mu_i \sim\mathcal{U}[-3,3]$, $i = 1,\ldots, 50$ & $\mathcal
{N}(\mu_i, 1)$ truncated on $[-5,5]$ \\
3 &
$\log(\sigma_i) \sim\mathcal{U}[-1,1]$, $\mu_i
\sim
\mathcal{U}[-2.5, 2.5]$, &$\mathcal{N}(\mu_i, \sigma_i^2)$ truncated on $[-5,5]$\\
& $\mu_i$ and $\sigma_i$ independent, $i = 1,
\ldots, 50$ &
\\
\hline
\end{tabular*}
\end{table}

In order to compare the different methods, we assessed the efficiency
of the resulting representations. Efficiency was quantified by the
fraction of variance explained (FVE), $\tilde{V}_K/\tilde{V}_\infty
$, as given by
the Fr\'echet variance [see (\ref{eq: tot_var_est}) and (\ref{eq:
var_ex_est})], so that higher FVE values reflect superior
representations. As this quantity depends on the chosen metric $d$, we
computed these values for both the $L^2$ and Wasserstein metrics. The
FVE results for the two metrics were similar, so we only present the
results using the $L^2$ metric here. Those corresponding to the
Wasserstein metric $d_W$ are given in the supplemental article \cite
{pete:15}. As mentioned in Section~\ref{sec: fda}, the truncated
representations in (\ref{eq: fnt_trnc}) given by ordinary FPCA are not
guaranteed to be bona fide densities. Hence, the representations were
first projected onto the space of densities by taking the positive part
and renormalizing, a method that has been systematically investigated
by \cite{gaje:86}.

\begin{figure}
\centering
\begin{tabular}{@{}ccc@{}}

\includegraphics{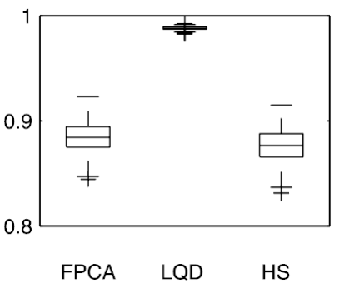}
 & \includegraphics{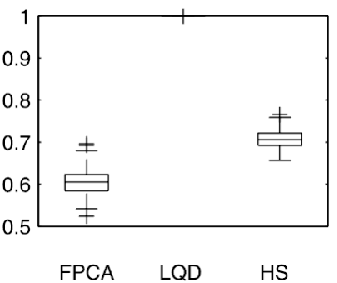}& \includegraphics{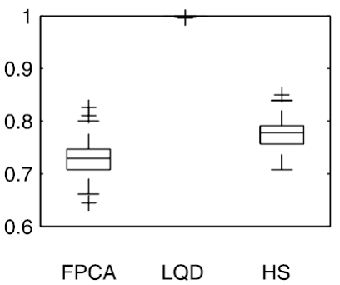}\\
\footnotesize{(a) Setting $1 - K = 1$} & \footnotesize{(b) Setting $2 - K =
1$}&\footnotesize{(c) Setting $3 - K = 2$}\\[3pt]

\includegraphics{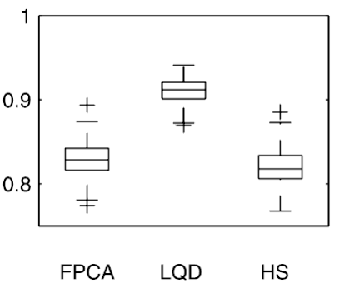}
& \includegraphics{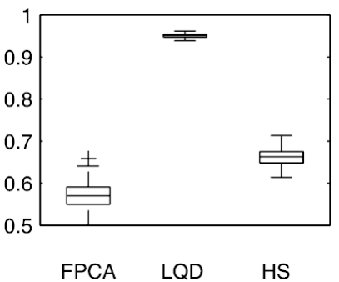} & \includegraphics{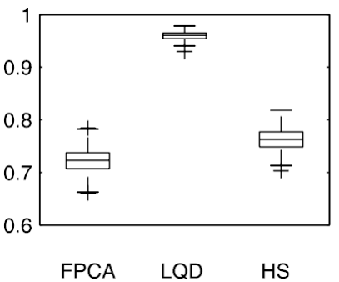}\\
\footnotesize{(d) Setting $1 - K = 1$}&
\footnotesize{(e) Setting $2 - K = 1$} & \footnotesize{(f) Setting $3 - K =
2$}
\end{tabular}
\caption{Boxplots of FVE (fraction of Fr\'echet variance
explained, larger is better) values for 200 simulations, using the
$L^2$ distance $d_2$. The first row corresponds to fully observed
densities and the second corresponds to estimated densities. The
columns correspond to settings 1, 2 and 3 from left to right (see
Table~\protect\ref{tab: sim_set}). The methods are denoted by
``FPCA'' for
ordinary FPCA on the densities, ``LQD'' for the transformation approach
with $\psi_Q$ and ``HS'' for the Hilbert sphere method.}
\label{fig: sim_fve}
\end{figure}

Boxplots for the FVE values (using the metric $d_2$) for the three
simulation settings are shown in Figure~\ref{fig: sim_fve}, where the
first row corresponds to fully observed densities and the second row to
estimated densities. The number of components used to compute the
fraction of variance explained was $K = 1$ for settings 1 and~2, and $K
= 2$ for setting 3, reflecting the true dimensions of the random
process generating the densities. Even in the first simulation setting,
where the variation is strictly vertical, the transformation method
outperformed both the standard FPCA and Hilbert sphere methods. The
advantage of the transformation is most noticeable in settings 2 and 3
where horizontal variation is prominent.

As a qualitative comparison, we also computed the Fr\'echet means
corresponding to three metrics: The $L^2$ metric (cross-sectional
mean), Wasserstein metric and Fisher--Rao metric. This last metric
corresponds to the geodesic metric on the Hilbert sphere between
square-root densities. This fact was exploited in \cite{sriv:07}, where
an estimation algorithm was introduced that we have implemented in our
analyses. For details on the estimation of the Wasserstein--Fr\'echet
mean, see the supplemental article \cite{pete:15}. To summarize these
mean estimates across simulations, we again took the Fr\'echet mean
(i.e., a Fr\'echet mean of Fr\'echet means), using the respective metric.

\begin{figure}
\centering
\begin{tabular}{@{}ccc@{}}

\includegraphics{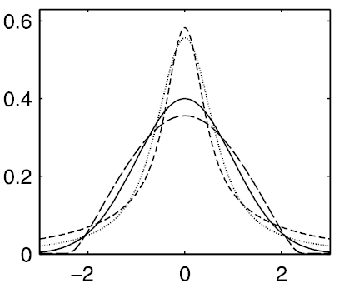}
 & \includegraphics{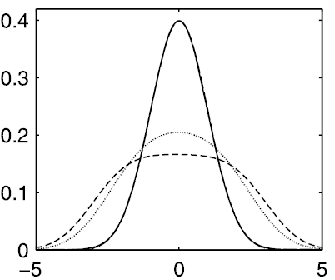}&\includegraphics{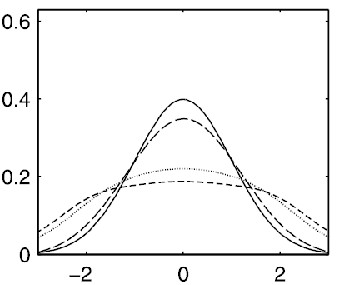} \\
\footnotesize{(a) Setting 1} & \footnotesize{(b) Setting 2}& \footnotesize
{(c) Setting 3}\\[3pt]

\includegraphics{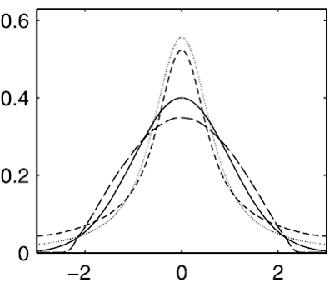}
&\includegraphics{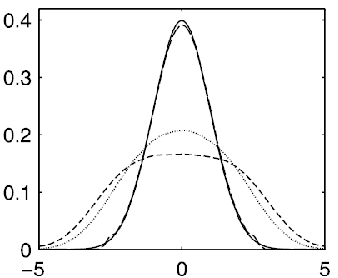} & \includegraphics{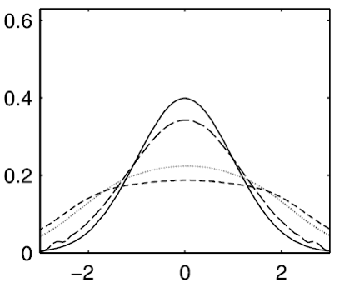}\\
\footnotesize{(d) Setting 1}& \footnotesize{(e) Setting 2} & \footnotesize
{(f) Setting 3}
\end{tabular}
\caption{Average Fr\'echet means across 200 simulations. The
first row corresponds to fully observed densities and the second
corresponds to estimated densities. The columns correspond to settings
1, 2 and 3 from left to right (see Table~\protect\ref{tab: sim_set}).
Truncated
$\mathcal{N}(0,1)$---solid line; Cross-sectional---short-dashed line;
Fisher--Rao---dotted line; Wasserstein---long-dashed line.}
\label{fig: sim_means}
\end{figure}

Note that a natural center for each simulation, if one knew the true
random mechanism generating the densities, is the (truncated) standard
normal density. Figure~\ref{fig: sim_means} plots the average mean
estimates across all simulations (in the Fr\'echet sense) for the
different settings along with the truncated standard normal density.
One finds that in setting 2 for fully observed densities, the
Wasserstein--Fr\'echet mean is visually indistinguishable from
truncated normal density. Overall, it is clear that the
Wasserstein--Fr\'echet mean yields a better concept for the ``center'' of
the distribution of data curves than either the cross-sectional or
Fisher--Rao--Fr\'echet means.

\subsection{Intra-hub connectivity and cognitive ability}
\label{ss: brain}

In recent years, the problem of identifying functional connectivity
between brain voxels or regions has received a great deal of attention,
especially for resting state fMRI \cite{alle:12,ferr:13,shel:13}.
Subjects are asked to relax while undergoing a fMRI brain scan, where
blood-oxygen-level dependent (BOLD) signals are recorded and then
processed to yield voxel-specific time courses of signal strength.
Functional connectivity between voxels is customarily quantified in
this area by the Pearson product-moment correlation \cite
{acha:06,bass:06,wors:05} which, from a functional data analysis point
of view,
corresponds to a special case of dynamic correlation for random
functions \cite{mull:05:2}. These correlations can be used for a variety
of purposes. A traditional focus has been on characterizing voxel
regions that have high correlations \cite{buck:09}, which have been
referred to as ``hubs.'' For each such hub, a so-called seed voxel is
identified as the voxel with the signal that has the highest
correlation with the signals of nearby voxels.

As a novel way to characterize hubs, we analyzed the distribution of
the correlations between the signal at the seed voxel of a hub and the
signals of all other voxels within an $11\times 11\times 11$ cube of voxels that is
centered at the seed voxel. For each subject, the target is the density
within a specified hub that is then estimated from the observed
correlations. The resulting sample of densities is then an i.i.d.
sample across subjects. To demonstrate our methods, we select the Right
inferior/superior Parietal Lobule hub (RPL) that is thought to be
involved in higher mental processing \cite{buck:09}.

The signals for each subject were recorded over the interval [0, 470]
(in seconds), with 236 measurements available at 2 second intervals.
For the fMRI data recorded for $n=68$ subjects that were diagnosed with
Alzheimer's disease at UC Davis, we performed standard preprocessing
that included the steps of slice-time correction, head motion
correction and normalization to the Montreal Neurological Institute
(MNI) fMRI template, in addition to linear detrending to account for
signal drift, band-pass filtering to include only frequencies between
0.01 and 0.08~Hz and regressing out certain time-dependent covariates
(head motion parameters, white matter and CSF signal).

\begin{figure}

\includegraphics{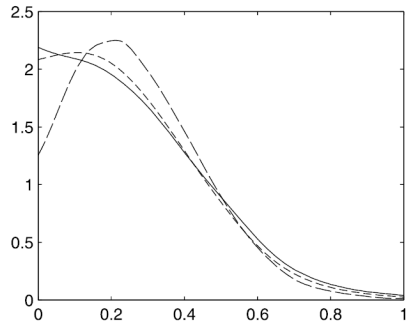}

\caption{Comparison of means for distributions of seed voxel
correlations for the RPL hub. Cross-sectional mean---solid line;
Fisher--Rao--Fr\'echet mean---short-dashed line; Wasserstein--Fr\'
echet mean---long-dashed line.}
\label{fig: rpl_means}
\end{figure}

For the estimation of the densities of seed voxel correlations, the
density estimator in (\ref{eq: dens_est}) was utilized, with kernel
$\kappa $ chosen as the standard Gaussian density and a bandwidth of
$h =
0.08$. As negative correlations are commonly ignored in connectivity
analyses, the densities were estimated on $[0,1]$. Figure~\ref{fig:
rpl_cor} shows the estimated densities for all $68$ subjects. A notable
feature is the variation in the location of the mode, as well as the
associated differences in the sharpness of the density at the mode. The
Fr\'echet means that one obtains with different approaches are plotted
in Figure~\ref{fig: rpl_means}. As in the simulations, the
cross-sectional and Fisher--Rao--Fr\'echet means are very similar, and
neither reflects the characteristics of the distributions in the
sample. In contrast, the Wasserstein--Fr\'echet mean displays a sharper
mode of the type that is seen in the sample of densities. Therefore, it
is clearly more representative of the sample.

\begin{figure}
\centering
\begin{tabular}{@{}ccc@{}}

\includegraphics{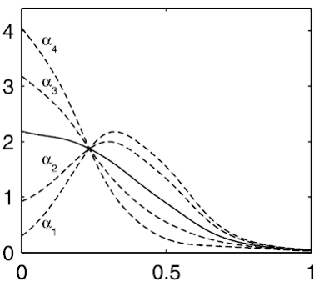}
 & \includegraphics{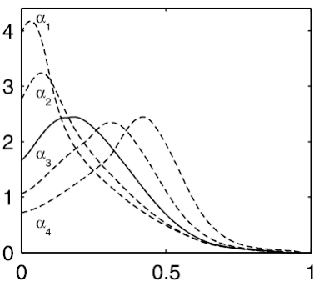}&\includegraphics{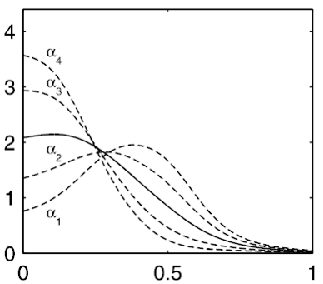} \\
\footnotesize{(a) Ordinary FPCA} & \footnotesize{(b) Log quantile density
transformation}&\footnotesize{(c) Hilbert sphere method} \\[3pt]

\includegraphics{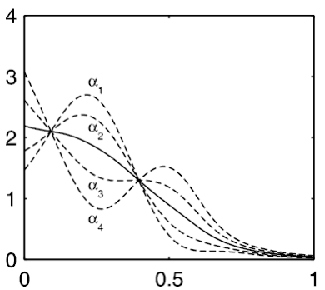}
& \includegraphics{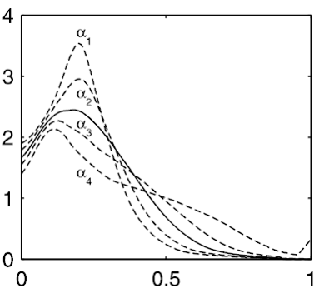} & \includegraphics{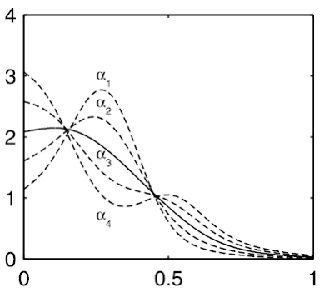}\\
\footnotesize{(d) Ordinary FPCA} & \footnotesize{(e) Log quantile density
transformation} & \footnotesize{(f) Hilbert sphere method}
\end{tabular}
\caption{Modes of variation for distributions of seed voxel
correlations. The first row corresponds to the first mode and the
second row to the second mode of variation. The values of $\alpha$ used
in the computation of the modes are quantiles ($\alpha_1 = 0.1$,
$\alpha
_2 = 0.25$, $\alpha_3 = 0.75$, $\alpha_4 = 0.9$) of the standardized
estimates of the principal component (geodesic) scores for each method,
and the solid line corresponds to $\alpha= 0$.}
\label{fig: rpl_modes}
\end{figure}

%
%

Next, we examined the first and second modes of variation, which are
shown in Figure~\ref{fig: rpl_modes}. The first mode of variation for
each method reflects the horizontal shifts in the density modes, the
location of which varies by subject. The modes for the Hilbert sphere
method closely resemble those for ordinary FPCA and both FPCA and
Hilbert sphere modes of variation do not adequately reflect the nature
of the main variability in the data, which is the shift in the modes
and associated shape changes. In contrast, the transformation modes of
variation using the log quantile density transformation retain the
sharp peaks seen in the sample and give a clear depiction of the
horizontal variation. The second mode describes vertical variation.
Here, the superiority of the transformation modes is even more
apparent. The modes of ordinary FPCA and, to a lesser extent, those for
the Hilbert sphere method, capture this form of variation awkwardly,
with the extreme values of $\alpha$ moving toward bimodality---a
feature that is not present in the data. In contrast, the log quantile
density modes of variation capture the variation in the peaks
adequately, representing all densities as unimodal density functions,
where unimodality is clearly present throughout the sample of density estimates.

\begin{figure}

\includegraphics{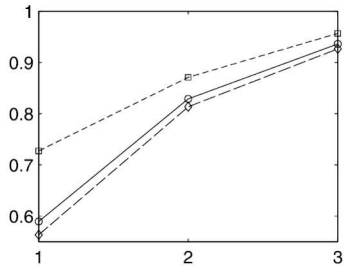}

\caption{Fraction of variance explained for $K = 1,2,3$
components, using the metric $d_2$. Ordinary FPCA---solid line/circle
marker; log quantile density transformation---short-dashed line/square
marker; Hilbert Sphere method---long-dashed line/diamond
marker.}
\label{fig: rpl_fve}
\end{figure}

In terms of connectivity, the first transformation mode reflects mainly
horizontal shifts in the densities of connectivity with associated
shape changes that are less prominent, and can be characterized as
moving from low to higher connectivity. The second transformation mode
of variation provides a measure of the peakedness of the density, and
thus to what extent connectivity is focused around a central value. The
fraction of variance explained as shown in Figure~\ref{fig: rpl_fve}
demonstrates that the transformation method provides not only more
interpretable modes of variation, but also more efficient
representations of the distributions than both ordinary FPCA and the
Hilbert sphere methods. Thus, while the transformation modes of
variation provide valuable insights into the variation of connectivity
across subjects, this is not the case for the ordinary or Hilbert
sphere modes of variation.

We also compared the utility of the densities and their transformed
versions to predict a cognitive test score which assesses executive
performance in the framework of a functional linear regression model.
As the Hilbert sphere method does not give a linear representation, it
cannot be used in this context. Denote the densities by $f_i$ with
functional principal components $\xi_{ik}$, the log quantile density
functions by $X_i = \psi_Q(f_i)$ with functional principal components
$\eta_{ik}$ and the test scores by $Y_i$. Then the two models \cite
{cai:06,hall:07:1} are
\begin{eqnarray*}
Y_i &=& B_{10} + \sum_{k = 1}^\infty
B_{1k}\xi_{ik} + \varepsilon _{1i} \quad\mbox{and}
\\
Y_i &=& B_{20} + \sum_{k = 1}^\infty
B_{2k}\eta_{ik} + \varepsilon _{2i},\qquad  i = 1,
\ldots,65,
\end{eqnarray*}
where three subjects who had missing test scores were removed. In
practice, the sums are truncated in order to produce a model fit. These
models were fit for different values of the truncation parameter $K$
[see (\ref{eq: fnt_trnc}) and (\ref{eq: fnt_trnc1})] using the PACE
package for MATLAB (code available at \url{http://anson.ucdavis.edu/\textasciitilde mueller/data/pace.html}) and 10-fold cross
validation (averaged over 50 runs) was used to obtain the mean squared
prediction error estimates give in Table~\ref{tab: rpl_mspe}.

\begin{table}[b]
\caption{Estimated mean squared prediction errors as obtained by
10-fold cross validation, averaged over 50 runs. Functional $R^2$
values for the fitted model using all data points are given in
parentheses}\label{tab: rpl_mspe}
\begin{tabular*}{\textwidth}{@{\extracolsep{\fill}}lcccc@{}}
\hline
$\bolds{K}$ & \textbf{1} & \textbf{2} & \textbf{3} & \textbf{4} \\
\hline
FPCA & 0.180 (0.0031) & 0.185 (0.0135) & 0.193 (0.0233) & 0.201 (0.0244) \\
LQD & 0.180 (0.0030) & 0.176 (0.0715) & 0.169 (0.1341) & 0.173 (0.1431) \\
\hline
\end{tabular*}
\end{table}

In addition, the models were fitted using all data points to obtain an
$R^2$ goodness-of-fit measurement for each truncation value $K$. The
transformed densities were found to be better predictors of executive
function than the ordinary densities for all values of $K$, both in
terms of prediction error and $R^2$ values. While the $R^2$ values were
generally small, as only a relatively small fraction of the variation
of the cognitive test score can generally be explained by connectivity,
they were much larger for the model that used the transformation scores
as predictors. These regression models relate transformation components
of brain connectivity to cognitive outcomes, and thus shed light on the
question of how patterns of intra-hub connectivity relate to cognitive function.

\section{Discussion}
\label{sec: discussion}

Due to the nonlinear nature of the space of density functions, ordinary
FPCA is problematic for functional data that correspond to densities,
both theoretically and practically, and the alternative transformation
methods as proposed in this paper are more appropriate. The
transformation based representations always satisfy the constraints of
the density space and retain a linear interpretation in a suitably
transformed space. The latter property is particularly useful for
functional regression models with densities as predictors. Notions of
mean and fraction of variance explained can be extended by the
corresponding Fr\'echet quantities once a metric has been chosen. The
Wasserstein metric is often highly suitable for the modeling of samples
of densities.

While it is well known that for the $L^2$ metric $d_2$ the
representations provided by ordinary FPCA are optimal in terms of
maximizing the fraction of explained variance among all $K$-dimensional
linear representations using orthonormal eigenfunctions, this is not
the case for other metrics or if the representations are constrained to
be in density space. In the transformation approach, the usual notion
of explained variance needs to be replaced. We propose to do this by
adopting the Fr\'echet variance, which in general will depend on the
chosen transformation space and metric. As the data analysis indicates,
even in the case of the $L^2$ metric, the log quantile density
transformation performs better compared to FPCA or the Hilbert sphere
approach in explaining most of the variation in a sample of densities
by the first few components. The FVE plots, as demonstrated in
Section~\ref{sec: application}, provide a convenient characterization
of the quality of a transformation and can be used to compare multiple
transformations or even to determine whether or not a transformation is
better than no transformation.

In terms of interpreting the variation of functional density data, the
transformation modes of variation emerge as clearly superior in
comparison to the ordinary modes of variation, which do not keep the
constraints to which density functions are subject. Overall, ordinary
FPCA emerges as ill-suited to represent samples of density functions.
When using such representations as an intermediate step, for example,
if prediction of an outcome or classification with densities as
predictors is of interest, it is likely that transformation methods are
often preferable, as demonstrated in our data example.

Various transformations can be used that satisfy certain continuity
conditions that imply consistency. In our experience, the log quantile
density transformation emerges as the most promising of these. While we
have only dealt with one-dimensional densities in this paper,
extensions to densities with more complex support are possible. Since
hazard and quantile functions are not immediately generalizable to
multivariate densities, there is no obvious extension of the
transformations based on these concepts to the multivariate case.
However, for multivariate densities, a relatively straightforward
approach is to apply the one-dimensional methodology to the conditional
densities used by the Rosenblatt transformation \cite{rose:52} to
represent higher-dimensional densities, although this approach would be
computationally demanding and is subject to the curse of dimensionality
and reduced rates of convergence as the dimension increases. However,
it would be quite feasible for two- or three-dimensional densities. In
general, the transformation approach is flexible, as it can be adopted
for any transformation that satisfies some regularity conditions and
maps densities to a Hilbert space.

\begin{appendix}\label{app}
\section*{Appendix: Details on theoretical results}
\subsection{Proofs of propositions and theorems}

This section contains proofs of Propositions~\ref{prop: dens_est} and
\ref{prop: transformations} and Theorems~\ref{thmm: modes} and \ref{thmm:
choice_K}. We also include some auxiliary lemmas. Additional proofs and
a complete listing of all assumptions can be found in~\cite{pete:15}.

\begin{pf*}{Proof of Proposition~\ref{prop: dens_est}}
Clearly, $\check{f}\geq0$ and $\int_0^1 \check{f}(x) \,dx = 1$. Set
\[
\mathring{f}(x) = \frac{1}{Nh}\sum_{l = 1}^N
\kappa \biggl(\frac{x -
W_l}{h} \biggr)w(x,h),
\]
so that $\check{f}= \mathring{f}/\int\mathring{f}$. Set $c_{\kappa
} =  (\int_0^1\kappa (u)
\,du )^{-1}$. For any $x\in[0,1]$ and $h < 1/2$, we have $1 \leq w(x,h)
\leq c_{\kappa } $, so that
\[
c_{\kappa }^{-1}\leq\inf_{y \in[0,1]} \int
_{-yh^{-1}}^{(1-y)h^{-1}
}\kappa (u) \,du\leq\int
_0^1 \mathring{f}(x) \,dx \leq c_{\kappa }.
\]
This implies
\[
\biggl\llvert 1 - \biggl(\int_0^1
\mathring{f}(x) \,dx \biggr)^{-1}\biggr\rrvert \leq\min \bigl\{
c_{\kappa } - 1, c_\kappa d_2(\mathring{f}, f),
c_\kappa d_{\infty
}(\mathring{f}, f) \bigr\},
\]
which,
together with assumption (A1), implies
\[
d_2(\check{f}, f) \leq c_\kappa (M + 1)\,d_2(
\mathring{f}, f) \quad\mbox {and}\quad d_{\infty} (\check{f}, f) \leq
c_\kappa (M + 1)d_{\infty}(\mathring{f}, f).
\]
Thus, we only need prove the remaining requirements in assumptions (D1)
and (D2) for the estimator $\mathring{f}$.

The expected value is given by
\begin{eqnarray*}
E\bigl(\mathring{f}(x)\bigr) &= &h^{-1}\int_0^1
\kappa \biggl(\frac{x -
y}{h} \biggr)w(x,h)f(y) \,dy
\\
&=& f(x) + h w(x,h)\int_{-xh^{-1}}^{(1-x)h^{-1}}f'
\bigl(x^\ast\bigr)u \kappa (u) \,dv,
\end{eqnarray*}
for some $x^\ast$ between $x$ and $x + uh$. Thus, $E(\mathring{f}(x))
= f(x) +
O(h)$, where the $O(h)$ term is uniform over $x \in[0,1]$ and $f\in
\mathcal
{F}$. Here, we have used the fact that $\sup_{f \in\mathcal {F}}
\Vert f' \Vert_\infty <
M$ and $\int_{\mathbb{R}} |u| \kappa (u) \,du < \infty$. Similarly,
\[
\operatorname{Var}\bigl(\mathring{f}(x)\bigr) \leq\frac{c_{\kappa
}^2}{Nh} \biggl(f(x)
\int_0^1\kappa ^2(u) \,du + h \int
_0^1 u\kappa ^2(u)f'
\bigl(x^\ast\bigr) \,du \biggr),
\]
for some $x^\ast$ between $x$ and $x + uh$, so that the variance is of
the order $(Nh)^{-1}$ uniformly over $x \in[0,1]$ and $f \in\mathcal {F}$.
This proves (D1) for $b_N^2 = h^2 + (Nh)^{-1}$.

To prove assumption (D2), we use the triangle inequality to see that
\[
d_{\infty}(f, \mathring{f}) \leq d_{\infty}\bigl(f, E\bigl(
\mathring{f}(\cdot )\bigr)\bigr) + d_{\infty}\bigl(\mathring{f}, E\bigl(
\mathring{f} (\cdot)\bigr)\bigr).
\]
Using the DKW inequality \cite{dvor:56}, there are constants $c_1$, $c_2$
and a sequence \mbox{$L_h = O(h)$} such that, for any $R > 0$,
\[
P\bigl(d_{\infty}(f, \mathring{f}) > 2Ra_N\bigr) \leq
c_1\exp \bigl\{ -c_2R^2a_N^2Nh^2
\bigr\} + I\{L_h > Ra_N\},
\]
where $I$ is the indicator function. Notice that the bound is
independent of $f \in\mathcal {F}$. By taking $h = N^{-1/3}$ and $a_N =
N^{-c}$ for $c \in(0, 1/6)$, we have $L_h < Ra_N$ for large enough
$N$, and thus, for such $N$,
\[
\sup_{f\in\mathcal {F}} P\bigl(d_{\infty}(f, \mathring{f}) >
2Ra_N\bigr) \leq c_1 \exp \bigl\{ -c_2R^2N^{1/3 - 2c}
\bigr\} = o(1) \qquad\mbox{as } N \rightarrow \infty.
\]
In assumption (S1), we may then take $m = n^r$ for any $r > 0$, since
\begin{eqnarray}
n\sup_{f\in\mathcal {F}} P\bigl(d_{\infty}(f, \mathring{f}) >
2Ra_N\bigr) \leq c_1 n\exp \bigl\{ -c_2R^2n^{r/3 - 2rc}
\bigr\} = o(1) \nonumber\\
\eqntext{\mbox{as } n \rightarrow \infty.\qquad}
\end{eqnarray}
\upqed\end{pf*}



\begin{pf*}{Proof of Proposition~\ref{prop: transformations}}
First, we deal with the log hazard transformation. Let $f$ and $g$ be
two densities as specified in assumption (T0), with distribution
functions $F$ and $G$. Then
\[
d_{\infty}(F, G) \leq d_2(f, g) \leq d_{\infty}(f, g).
\]
Also, $1 - F$ and $1 - G$ are both bounded below by $\delta D_0^{-1}$ on
$[0, {1_\delta}]$. Then, for $x \in[0,{1_\delta}]$,
\begin{eqnarray*}
\bigl\llvert \psi_H(f) (x) - \psi_H(g) (x)\bigr\rrvert
&\leq&\biggl\llvert \log \biggl(\frac
{f(x)}{g(x)} \biggr)\biggr\rrvert + \biggl
\llvert \log \biggl(\frac{1 - F(x)}{1 -
G(x)} \biggr)\biggr\rrvert
\\
&\leq& D_0 \bigl[\bigl\llvert f(x) - g(x)\bigr\rrvert +
\delta^{-1}\bigl\llvert F(x) - G(x)\bigr\rrvert \bigr],
\end{eqnarray*}
whence
\begin{eqnarray*}
d_{\infty}\bigl(\psi_H(f), \psi_H(g)\bigr) &\leq&
D_0\bigl(1 + \delta ^{-1}\bigr)d_{\infty}(f, g),
\\
d_2\bigl(\psi_H(f), \psi_H(g)
\bigr)^2 &\leq&2D_0^2 \biggl[\int
_0^{{1_\delta}
}\bigl(f(x) - g(x)\bigr)^2 \,dx +
\delta^{-2}d_2(f, g)^2 \biggr]
\\
&\leq&2D_0^2 \bigl(1 + \delta^{-2}
\bigr)d_2(f, g)^2.
\end{eqnarray*}
These bounds provide the existence of $C_0$ in (T0). For (T1), observe
that
\[
\delta D_1^{-2} < \frac{f(x)}{1 - F(x)} \leq
\delta^{-1}D_1^2,
\]
so that
\begin{eqnarray*}
\bigl\Vert\psi_H(f) \bigr\Vert_\infty &=& \sup_{x \in[0,{1_\delta}]}
\biggl\llvert \log\frac
{f(x)}{1 -
F(x)}\biggr\rrvert \leq2\log D_1 - \log
\delta\quad \mbox{and}
\\
\bigl\Vert\psi_H(f)' \bigr\Vert_\infty &=& \sup
_{x \in[0, {1_\delta}]}\biggl\llvert \frac{f'(x)(1
- F(x)) +
f(x)^2}{f(x)(1 - F(x))}\biggr\rrvert \leq2
\delta^{-1}D_1^4,
\end{eqnarray*}
which proves the existence of $C_1$.

Next, let $X$ and $Y$ be functions as in (T2) for $\mathcal {T} = [0,
{1_\delta}]$
and set \mbox{$f = \psi_H^{-1}(X)$} and \mbox{$g = \psi_H^{-1}
(Y)$}. Let
$\Lambda_X(x) = \int_0^x e^{X(s)} \,ds$ and $\Lambda_Y(x) = \int_0^xe^{Y(s)} \,ds$. Then
\[
\bigl\llvert \Lambda_X(x) - \Lambda_Y(x)\bigr\rrvert
\leq\int_0^x\bigl\llvert e^{X(s)} -
e^{Y(s)}\bigr\rrvert \,ds \leq e^{\Vert X \Vert_\infty + d_{\infty}(X, Y)}d_2(X,Y),
\]
whence
%
\begin{eqnarray}\label{eq: l2lh_bound}
&&d_2\bigl(\psi_H^{-1}(X),
\psi_H^{-1}(Y)\bigr)^2 \nonumber\\
&&\qquad\leq2e^{2\Vert X \Vert
_\infty}
\bigl[d_2 (\Lambda_X, \Lambda_Y
)^2 \,dx + e^{2d_{\infty}(X, Y)}d_2(X, Y)^2 \bigr]
\nonumber
\\[-8pt]
\\[-8pt]
\nonumber
&&\quad\qquad{} + \delta^{-1} \bigl(\Lambda_X({1_\delta}) -
\Lambda _Y({1_\delta}) \bigr)^2
\nonumber
\\
&&\qquad\leq2e^{2\Vert X \Vert_\infty} \bigl[ \bigl(e^{2\Vert X \Vert
_\infty} + \delta^{-1}
\bigr) + 1 \bigr]e^{2d_{\infty}(X, Y)}d_2(X, Y)^2
.\nonumber
\end{eqnarray}
Taking \mbox{$C_2 = \sqrt{2}e^{\Vert X \Vert_\infty} [
(e^{2\Vert X \Vert_\infty} +
\delta^{-1} ) + 1 ]^{1/2}$} and \mbox{$C_3 = e^{d_{\infty
}(X, Y)}$},
(T2) is established for $d = d_2$.

For $d = d_W$, the cdf's of $f$ and $g$ for $x \in[0, {1_\delta}]$ are given
by $F(x) = 1 - e^{-\Lambda_X(x)}$ and \mbox{$G(x) = 1 -
e^{-\Lambda_Y(x)}$}, respectively. For $x \in({1_\delta}, 1]$,
\begin{eqnarray*}
F(x) &=& F({1_\delta}) + \delta^{-1}\bigl(1 -
F({1_\delta})\bigr) (x - {1_\delta }),\\
 G(x)& =&
G({1_\delta}) + \delta^{-1}\bigl(1 - G({1_\delta})
\bigr) (x - {1_\delta}),
\end{eqnarray*}
so that $|F(x) - G(x)| \leq|F({1_\delta}) - G({1_\delta})|$ for such $x$.
Hence, for
all $x \in[0, 1]$
\[
\bigl\llvert F(x) - G(x)\bigr\rrvert \leq\sup_{x \in[0, {1_\delta}]}\bigl
\llvert \Lambda_X(x) - \Lambda_Y(x)\bigr\rrvert \leq
e^{\Vert X \Vert_\infty + d_{\infty}(X,
Y)}d_2(X,Y).
\]
Note that for $t \in[0, 1]$ and $t \neq F({1_\delta})$,
\[
\bigl(F^{-1}\bigr) '(t) = \bigl[f\bigl(F^{-1}(t)
\bigr) \bigr]^{-1}\leq\exp \bigl\{ e^{\Vert X \Vert_\infty} \bigr\}\max \bigl(
\delta^{-1}, e^{\Vert X
\Vert_\infty} \bigr)=: c_L,
\]
so that $F^{-1}$ is Lipschitz with constant $c_L$. Thus, letting $t \in
[0,1]$ and $x = G^{-1}(t)$,
\[
\bigl\llvert F^{-1}(t) - G^{-1}(t)\bigr\rrvert = \bigl
\llvert F^{-1}\bigl(G(x)\bigr) - F^{-1} \bigl(F(x)\bigr)\bigr
\rrvert \leq c_Le^{\Vert X \Vert_\infty + d_{\infty}(X, Y)}d_2(X,Y),
\]
whence
%
\begin{equation}\qquad
\label{eq: wlh_bound} d_W\bigl(\psi_H^{-1}(X),
\psi_H^{-1}(Y)\bigr) = d_2\bigl(F^{-1},
G^{-1}\bigr) \leq c_Le^{\Vert X \Vert_\infty}e^{d_{\infty}(X, Y)}d_2(X,
Y).
\end{equation}
Using (\ref{eq: wlh_bound}),\vspace*{1pt} we establish (T2) for $d_W$ by setting
$C_2 = c_Le^{\Vert X \Vert_\infty}$ and $C_3 = e^{d_{\infty}(X, Y)}$.

To establish (T3), we let $X = \psi_H(f_1)$ and $X_K = \nu+ \sum_{k =
1}^K \eta_{1k} \rho_k$.
Set $f_{1,K} = \psi_H^{-1}(X_K)$ and take $C_1$ as in (T1). Then, by
assumption (A1) and equations (\ref{eq: l2lh_bound}) and (\ref{eq:
wlh_bound}),
\begin{eqnarray*}
E \bigl(d_2(f_1, f_{1,K})^2 \bigr)
&\leq& b_1\sqrt{E \bigl(e^{4d_{\infty}(X,
X_K)} \bigr)E \bigl(d_2(X,
X_K)^4 \bigr)} \quad \mbox{and}
\\
E\bigl(d_W(f_1, f_{1,K})^2\bigr)
&\leq &b_2\sqrt{E \bigl(e^{4d_{\infty}(X,
X_K)} \bigr)E \bigl(d_2(X,
X_K)^4 \bigr)},
\end{eqnarray*}
where $b_1 = 2e^{2C_1} [ (e^{2C_1} + \delta^{-1} ) +
1
]$ and $b_2 = \exp \{2(e^{C_1}+C_1) \}\max (\delta
^{-2},e^{2C_1} )$. Note that $d_2(X, X_K)^2 = \sum_{k =
K+1}^\infty\eta_{1k}^2 \leq\Vert X \Vert_2^2 \leq C_1^2|\mathcal
{T}|$, so that
\[
E\bigl(d_2(X, X_K)^4\bigr) \leq
C_1^2|\mathcal {T}|E \Biggl(\sum
_{k =
K+1}^\infty \eta _{1k}^2
\Biggr) = C_1^2|\mathcal {T}|\sum
_{k = K+1}^\infty\tau_k \rightarrow 0.
\]
So, we just need to show that $E (e^{4d_{\infty}(X, X_K)} )
= O(1)$.

For the following, we need two lemmas that are listed below, and whose
proofs are in the online supplement \cite{pete:15}.
By applying assumptions (A1) and (T1), Lemma~\ref{lma: fpca_bounds}
implies the existence of the Lipschitz constant $L_K$ for the residual
process $X - X_K$ [see (\ref{eq: lip_const})]. By Lemma~\ref{lma:
lipschitz}, we have
\[
E \bigl(e^{4d_{\infty}(X, X_K)} \bigr) \leq E \bigl(\exp \bigl\{ 8|A|^{-1/2}
d_2(X, X_K) \bigr\} + \exp \bigl\{8L_K^{1/3}d_2(X,X_K)^{2/3}
\bigr\} \bigr).
\]
Since $d_2(X, X_K) \leq\Vert X \Vert_2 < C_1|\mathcal {T}|^{1/2}$,
the first
expectation is bounded. For the second, we use Jensen's inequality to find
%
\begin{eqnarray}
\label{eq: step1} &&E \bigl(\exp \bigl\{8L_K^{1/3}d_2(X,X_K)^{2/3}
\bigr\} \bigr)
\nonumber
\\[-8pt]
\\[-8pt]
\nonumber
&&\qquad\leq 1 + \sum_{m = 1}^\infty
\frac{8^m  [L_K^m E(d_2(X,
X_K)^{2m}) ]^{1/3}}{m!}.
\end{eqnarray}
For r.v.s. $Y_1, \ldots, Y_m$, $E(\prod_{i = 1}^m Y_i) \leq\prod_{i =
1}^m E(Y_i^m)^{1/m}$, so that
\begin{eqnarray*}
E\bigl(d_2(X, X_K)^{2m}\bigr) &=& \sum
_{k_1 = K+1}^\infty\cdots\sum_{k_m =
K+1}^\infty
E \Biggl(\prod_{i = 1}^m
\eta_{1k_i}^2 \Biggr)
\\
&\leq&\sum_{k_1 = K+1}^\infty\cdots\sum
_{k_m = K+1}^\infty\prod_{i
= 1}^m
E\bigl(\eta_{1k_i}^{2m}\bigr)^{1/m} = \Biggl(\sum
_{k = K + 1}^\infty E\bigl(\eta _{1k}^{2m}
\bigr)^{1/m} \Biggr)^m.
\end{eqnarray*}
Next, by assumption, there exists $B$ such that $L_K\sum_{k =
K+1}^\infty\tau_k \leq B$ for large $K$. Then, by the assumption on
the higher moments of $\eta_{1k}^{2m}$, for large $K$
\begin{eqnarray*}
L_K^mE\bigl(d_2(X, X_K)^{2m}
\bigr) &\leq& \Biggl(L_K\sum_{k = K+1}^\infty
E\bigl(\eta _{1k}^{2m}\bigr)^{1/m}
\Biggr)^m \leq \Biggl(L_K\sum
_{k = K+1}^\infty\bigl(r_m
\tau_k^m\bigr)^{1/m} \Biggr)^m
\\
&\leq& r_mB^m.
\end{eqnarray*}
Inserting this into (\ref{eq: step1}), for large $K$
\[
E \bigl(\exp \bigl\{8L_K^{1/3}d_2(X,X_K)^{2/3}
\bigr\} \bigr) \leq 1 + \sum_{m = 1}^\infty
\frac{8^m B^{m/3}r_m^{1/3}}{m!}.
\]
Using the assumption that $ (\frac{r_{m+1}}{r_m} )^{1/3} =
o(m)$, the ratio test shows the sum converges. Since the sum is
independent of $K$ for $K$ large, this establishes that $E(d_W(f_1,
f_{1,K})^2) = o(1)$ and \mbox{$E(d_2(f_1, f_{1,K})^2) = o(1)$}. Using
similar arguments, we can show that $E(d_W(f_1, f_{1,K})^4)$ and
$E(d_2(f_1, f_{1,K})^4)$ are both $O(1)$, which completes the proof.

Next, we prove (T0)--(T3) for the log quantile density transformation.
Let $f$ and $g$ be two densities as specified in assumption (T0) with
cdf's $F$ and $G$. For $t \in[0,1]$,
\begin{eqnarray*}
&&\bigl\llvert \psi_Q(f) (t) - \psi_Q(g) (t)\bigr\rrvert
\\
&&\qquad= \bigl\llvert \log f\bigl(F^{-1} (t)\bigr) - \log g
\bigl(G^{-1}(t)\bigr)\bigr\rrvert
\\
&&\qquad\leq D_0 \bigl(\bigl\llvert f\bigl(F^{-1}(t)\bigr) - f
\bigl(G^{-1}(t)\bigr)\bigr\rrvert + \bigl\llvert f\bigl(G^{-1}
(t)\bigr) - g\bigl(G^{-1}(t)\bigr)\bigr\rrvert \bigr)
\\
&&\qquad\leq D_0^2\bigl\llvert F^{-1}(t) -
G^{-1}(t)\bigr\rrvert + D_0\bigl\llvert f
\bigl(G^{-1} (t)\bigr) - g\bigl(G^{-1}(t)\bigr)\bigr\rrvert .
\end{eqnarray*}
Since $F' = f$ is bounded below by $D_0^{-1}$, for any $t \in[0,1]$ and
$x = G^{-1}(t)$,
\[
\bigl\llvert F^{-1}(t) - G^{-1}(t)\bigr\rrvert = \bigl
\llvert F^{-1}\bigl(G(x)\bigr) - F^{-1} \bigl(F(x)\bigr)\bigr
\rrvert \leq D_0\bigl\llvert F(x) - G(x)\bigr\rrvert .
\]
Recall that \mbox{$d_{\infty}(F, G) \leq d_2(f, g) \leq d_{\infty
}(f, g)$}. Hence,
\begin{eqnarray*}
d_{\infty}\bigl(\psi_Q(f), \psi_Q(g)\bigr) &\leq&
D_0 \bigl(D_0^2 + 1 \bigr)d_{\infty}
(f, g),
\\
d_2\bigl(\psi_Q(f), \psi_Q(g)
\bigr)^2 &\leq&2D_0^2 \biggl[D_0^4d_2(f,
g)^2 + \int_0^1\bigl(f(x) - g(x)
\bigr)^2 g(x) \,dx \biggr]
\\
&\leq&2D_0^3\bigl(D_0^3 + 1
\bigr)d_2(f, g)^2,
\end{eqnarray*}
whence $C_0$ in (T0). Next, we find that
\[
\bigl\Vert\psi_Q(f) \bigr\Vert_\infty \leq\log D_1\quad
\mbox{and}\quad \bigl\Vert\psi _Q(f)' \bigr\Vert_\infty\leq
D_1^3,
\]
whence $C_1$ in (T1).

Now, let $X$ and $Y$ be as stated in (T2). Let $F$ and $G$ be the
quantile functions corresponding to \mbox{$f = \psi_Q^{-1}(X)$} and
$g =
\psi_Q^{-1}(Y)$, respectively. Then
\begin{eqnarray*}
\bigl\llvert F^{-1}(t) - G^{-1}(t)\bigr\rrvert &\leq&
\theta_X^{-1}\biggl\llvert \int_0^t
\bigl(e^{X(s)} - e^{Y(s)} \bigr) \,ds\biggr\rrvert + \bigl\llvert
\theta_X^{-1}- \theta _Y^{-1}\bigr
\rrvert \int_0^te^{Y(s)} \,ds
\\
&\leq&2\theta_X^{-1}\llvert \theta_X -
\theta_Y\rrvert ,
\end{eqnarray*}
where $\theta_X = \int_0^1 e^{X(s)} \,ds$ and $\theta_Y = \int_0^1e^{Y(s)}
\,ds$. It is clear that $\theta_X^{-1}\leq e^{\Vert X \Vert_\infty}$
and \mbox
{$|\theta
_X - \theta_Y| \leq e^{\Vert X \Vert_\infty + d_{\infty}(X,
Y)}d_2(X,Y)$}, whence
\[
\bigl\llvert F^{-1}(t) - G^{-1}(t)\bigr\rrvert
\leq2e^{2\Vert X \Vert_\infty +
d_{\infty}(X,
Y)}d_2(X, Y).
\]
This implies
%
\begin{equation}
\label{eq: wlq_bound} d_W\bigl(\psi_Q^{-1}(X),
\psi_Q^{-1}(Y)\bigr) \leq2e^{4\Vert X \Vert_\infty
}e^{2d_{\infty}
(X, Y)}d_2(X,Y).
\end{equation}
For $d = d_2$, using similar arguments as above, we find that
%
\begin{eqnarray}\label{eq: l2lq_bound}
&&d_2\bigl(\psi_Q^{-1}(X),
\psi_Q^{-1}(Y)\bigr)
\nonumber
\\[-8pt]
\\[-8pt]
\nonumber
&&\qquad \leq\sqrt{2}e^{6\Vert X \Vert
_\infty} \bigl(4
\bigl\Vert X' \bigr\Vert_\infty^2 + 3
\bigr)^{1/2}e^{2d_{\infty}(X, Y)}d_2(X, Y).
\end{eqnarray}
Equations (\ref{eq: wlq_bound}) and (\ref{eq: l2lq_bound}) can then be
used to find the constants $C_2$ and $C_3$ in (T2) for both $d = d_W$
and $d = d_2$, and also to prove (T3) in a similar manner to the log
hazard transformation.
\end{pf*}

The following auxiliary results, which are proved in the online
supplement, are needed.

\begin{lma}
\label{lma: lipschitz}
Let $A$ be a closed and bounded interval of length $|A|$ and assume
$X: A \rightarrow \mathbb {R}$ is continuous with Lipschitz constant
$L$. Then
\[
\Vert X \Vert_\infty \leq2\max \bigl(|A|^{-1/2}\Vert X
\Vert_2, L^{1/3}\Vert X \Vert_2^{2/3}
\bigr).
\]
\end{lma}

\begin{lma}
\label{lma: fpca_bounds}
Let $X$ be a stochastic process on a closed interval $\mathcal
{T}\subset
\mathbb{R}$ such that $\Vert X \Vert_\infty < C$ and \mbox{$\Vert
X' \Vert_\infty < C$} almost
surely. Let $\nu$ and $H$ be the mean and covariance functions
associated with $X$, and $\rho_k$ and $\tau_k$, $k\geq1$, be the
eigenfunctions and eigenvalues of the integral operator with kernel
$H$. Then $\Vert\nu \Vert_\infty < C$, $\Vert H \Vert_\infty <
4C^2$ and $\Vert\rho_k \Vert_\infty <
4C^2|\mathcal {T}|^{1/2}\tau_k^{-1}$ for all $k \geq1$. Additionally,
$\Vert\nu' \Vert_\infty < C$ and \mbox{$\Vert\rho_k' \Vert
_\infty < 4C^2|\mathcal {T}|^{1/2}\tau
_k^{-1}$}
for all $k \geq1$.
\end{lma}

\begin{lma}
\label{lma: mean_cov}
Under assumptions \textup{(A1)} and \textup{(T1)}, with $\hat{\nu}, \tilde{\nu},
\widehat {H}, \widetilde {H}$ as in
(\ref{eq: nu_est}) and~(\ref{eq: H_est}),
\begin{eqnarray*}
d_2(\nu, \tilde{\nu}) &=& O_p\bigl(n^{-1/2}
\bigr),\qquad  d_2(H, \widetilde {H}) = O_p\bigl(n^{-1/2}
\bigr),
\\
d_{\infty}(\nu, \tilde{\nu}) &= &O_p \biggl( \biggl(
\frac{\log
n}{n} \biggr)^{1/2} \biggr),\qquad  d_{\infty}(H, \widetilde
{H}) = O_p \biggl( \biggl(\frac{\log n}{n} \biggr)^{1/2}
\biggr).
\end{eqnarray*}
Under the additional assumptions \textup{(D1), (D2)} and \textup{(S1)}, we have
\begin{eqnarray*}
d_2(\nu, \hat{\nu}) &=& O_p\bigl(n^{-1/2}+
b_{m}\bigr),\qquad d_2(H, \widehat {H}) = O_p
\bigl(n^{-1/2}+ b_{m}\bigr),
\\
d_{\infty}(\nu, \hat{\nu}) &=& O_p \biggl( \biggl(
\frac{\log
n}{n} \biggr)^{1/2} + a_{m} \biggr),\qquad
d_{\infty}(H, \widehat {H}) = O_p \biggl( \biggl(
\frac{\log
n}{n} \biggr)^{1/2} + a_{m} \biggr).
\end{eqnarray*}
\end{lma}

\begin{lma}
\label{lma: fnt_trnc}
Assume \textup{(A1)}, \textup{(T1)} and \textup{(T2)} hold. Let $A_k = \Vert\rho_k \Vert_\infty
$, $M$ as in
\textup{(A1)}, $\delta_k$ as in (\ref{eq: sp}), and $C_1$ as in \textup{(T1)} with $D_1 =
M$. Let $K^\ast(n) \rightarrow \infty$ be any sequence which
satisfies $\tau
_{K^\ast}n^{1/2}\rightarrow \infty$ and
\[
\sum_{k = 1}^{K^\ast} \bigl[ (\log
n)^{1/2} + \delta_k^{-1}+ A_k + \tau
_{K^{\ast}}\delta_k^{-1}A_k \bigr] = O
\bigl(\tau_{K^{\ast}}n^{1/2}\bigr).
\]
Let $C_2$ be as in \textup{(T2)}, $X_{i,K}= \nu+ \sum_{k = 1}^K \eta_{ik}\rho_k$,
$\widetilde {X}_{i,K}= \tilde{\nu}+ \sum_{k = 1}^K \tilde{\eta
}_{ik}\tilde{\rho}_k$, and set
\[
S_{K^\ast} = \max_{1 \leq K \leq K^\ast}\max_{1 \leq i \leq n}
C_2\bigl(\Vert X_{i,K} \Vert_\infty, \bigl\Vert
X_{i,K}' \bigr\Vert_\infty\bigr).
\]
Then
\[
\max_{1 \leq K \leq K^\ast}\max_{1 \leq i \leq n}d\bigl(f_i(
\cdot, K, \psi), \tilde{f}_i(\cdot, K, \psi)\bigr) = O_p
\biggl(\frac{S_{K^\ast}\sum_{k =
1}^{K^\ast}
\delta_k^{-1}}{n^{1/2}} \biggr).
\]
\end{lma}

We now can also state the following corollary, the proof of which
utilizes a lemma from \cite{mull:08:2}.

\begin{cor}
\label{cor: fpca_est}
Under assumption \textup{(A1)} and \textup{(T1)}, letting $A_k = \Vert\rho_k \Vert
_\infty$, with
$\delta_k$ as in (\ref{eq: sp}),
\begin{eqnarray*}
|\tau_k - \tilde{\tau}_k| &=& O_p
\bigl(n^{-1/2}\bigr),
\\
d_2(\rho_k, \tilde{\rho}_k) &=&
\delta_k^{-1}O_p\bigl(n^{-1/2}\bigr)\quad
\mbox{and}
\\
d_{\infty}(\rho_k, \tilde{\rho}_k) &=& \tilde{
\tau}_k^{-1}O_p \biggl(\frac{(\log
n)^{1/2} +
\delta_k^{-1}+ A_k}{n^{1/2}}
\biggr),
\end{eqnarray*}
where all $O_p$ terms are uniform over $k$. If the additional
assumptions \textup{(D1)}, \textup{(D2)} and \textup{(S1)} hold,
\begin{eqnarray*}
|\tau_k - \hat{\tau}_k| &=& O_p
\bigl(n^{-1/2}+ b_{m}\bigr),
\\
d_2(\rho_k, \hat{\rho}_k) &=&
\delta_k^{-1}O_p\bigl(n^{-1/2}+
b_{m}\bigr)\qquad \mbox{and}
\\
d_{\infty}(\rho_k, \hat{\rho}_k) &=& \hat{
\tau}_k^{-1}O_p \biggl(\frac{(\log
n)^{1/2} +
\delta_k^{-1}+ A_k}{n^{1/2}} +
a_{m} + b_{m}\bigl[\delta_k^{-1}+
A_k\bigr] \biggr),
\end{eqnarray*}
where again all $O_p$ terms are uniform over $k$.
\end{cor}


\begin{pf*}{Proof of Theorem~\ref{thmm: modes}}
We will show the result for the fully observed case. The same arguments
apply to the case where the densities are estimated.

First, suppose $K$ is fixed. We may use the results of Lemma~\ref{lma:
fpca_bounds} due to (A1) and (T1) and define $A_k$ as in Corollary~\ref
{cor: fpca_est}. From
\[
Y_{k, \alpha} = \nu+ \alpha\sqrt{\tau_k}\rho_k\quad
\mbox{and}\quad \widetilde {Y}_{k, \alpha} = \tilde{\nu}+ \alpha\sqrt{\tilde {
\tau}_k}\tilde{\rho}_k,
\]
$g_k(\cdot, \alpha, \psi) = \psi^{-1}(Y_{k, \alpha})$ and
similarly for
$\tilde{g}_k$. Observe that, if $|\alpha| \leq\alpha_0$,
%
\begin{equation}
\label{eq: dinf_bound} d_{\infty}(Y_{k, \alpha}, \widetilde
{Y}_{k, \alpha}) \leq d_{\infty}(\nu, \tilde{\nu}) + \alpha
_0 \bigl(\sqrt{\tilde{\tau}_1}d_{\infty}(
\rho_k, \tilde{\rho}_k) + A_k\llvert \sqrt
{\tau_k} - \sqrt{\tilde{\tau}_k}\rrvert \bigr).
\end{equation}
Next, $\max_{1 \leq k \leq K}|\sqrt{\tau_k} - \sqrt{\tilde{\tau
}_k}| =
O_p(n^{-1/2})$ and \mbox{$\max_{1 \leq k \leq K}d_{\infty}(\rho_k,
\tilde{\rho}
_k) =
O_p(1)$} by Corollary~\ref{cor: fpca_est}, so that $d_{\infty}(Y_{k,
\alpha
}, \widetilde {Y}_{k, \alpha}) = O_p(1)$, uniformly in $k$ and
$|\alpha| \leq
\alpha_0$. For $C_{2,k, \alpha} = C_2(\Vert Y_{k, \alpha} \Vert
_\infty, \Vert Y_{k, \alpha}' \Vert_\infty)$ and $C_{3, k, \alpha}
= C_3(d_{\infty}(Y_{k, \alpha}, \widetilde {Y}_{k,
\alpha}))$ as in (T2),
\[
\max_{1 \leq k \leq K}\max_{|\alpha| \leq\alpha_0} C_{2, k, \alpha
} <
\infty\quad \mbox{and}\quad \max_{1 \leq k \leq K}\max_{|\alpha| \leq
\alpha_0}
C_{3, k, \alpha} = O_p(1).
\]
Furthermore,
\[
d_2(Y_{k, \alpha}, \widetilde {Y}_{k, \alpha}) \leq
d_2(\nu, \tilde {\nu}) + \alpha_0 \bigl(\sqrt{\tilde{
\tau}_1}d_2(\rho_k, \tilde{
\rho}_k) + \llvert \sqrt {\tau_k} - \sqrt{\tilde{\tau}
_k}\rrvert \bigr) = O_p\bigl(n^{-1/2}\bigr),
\]
uniformly in $k$ and $|\alpha| \leq\alpha_0$, by Lemma~\ref{lma:
mean_cov}. This means
\begin{eqnarray*}
\max_{1\leq k \leq K}\max_{|\alpha|\leq\alpha_0}d\bigl(g_k(
\cdot, \alpha, \psi), \tilde{g}_k(\cdot, \alpha, \psi)\bigr) &\leq&
\max_{1\leq k \leq
K}\max_{|\alpha|\leq\alpha_0} C_{2,k, \alpha}
C_{3, k, \alpha} d_2(Y_{k,\alpha}, \widetilde {Y}_{k,\alpha})
\\
&=& O_p\bigl(n^{-1/2}\bigr).
\end{eqnarray*}

Next, we consider $K = K(n) \rightarrow \infty$. Define
\[
S_{K} = \max_{|\alpha| \leq\alpha_0}\max_{1 \leq k \leq K}
C_{2, k,
\alpha}.
\]
Let $B_K = \max_{1 \leq k \leq K} A_k$ and take $K$ to be a sequence
which satisfies:
\begin{enumerate}[(iii)]
\item[(i)]$\tau_{K}n^{1/2}\rightarrow \infty$,
\item[(ii)]$(\log n)^{1/2}+ \delta_{K}^{-1}+ B_{K} = O(\tau_{K} n^{1/2})$, and
\item[(iii)]$S_{K} = o (\delta_{K} n^{1/2} )$.
\end{enumerate}
For $|\alpha|\leq\alpha_0$, we still have inequality (\ref{eq:
dinf_bound}). The term $d_{\infty}(\nu, \tilde{\nu})$ is $o_p(1)$
independently of
$K$. From (i) and the above, it follows that \mbox{$\max_{1 \leq k
\leq
K} \tilde{\tau}_k^{-1}= O_p(\tau_{K}^{-1})$} and we find
\[
\max_{1 \leq k \leq K}\llvert \sqrt{\tau_k} - \sqrt{\tilde{
\tau }_k}\rrvert = O_p \biggl(\frac{1}{(\tau_{K} n)^{1/2}}
\biggr).
\]
Using Corollary~\ref{cor: fpca_est} and (ii), this implies $\max_{1
\leq k \leq K} d_{\infty}(\rho_k, \tilde{\rho}_k) = o_p(1)$, so
that \mbox
{$d_{\infty}
(Y_{k, \alpha}, \widetilde {Y}_{k, \alpha}) = O_p(1)$}, uniformly
over $k \leq K$
and $|\alpha| \leq\alpha_0$. Hence,  $\max_{1 \leq k \leq K}\max_{|\alpha| \leq\alpha_0} C_{3, k, \alpha} = O_p(1)$.

Similarly, we find that
\[
d_2(Y_{k, \alpha}, \widetilde {Y}_{k, \alpha}) =
O_p \biggl(\frac
{1}{\delta
_{K}n^{1/2}
} \biggr),
\]
uniformly over $k \leq K(n)$ and $|\alpha| \leq\alpha_0$. With (iii),
this yields
\[
\max_{|\alpha| \leq\alpha_0}\max_{1 \leq k \leq K}d\bigl(g_k(
\cdot, \alpha, \psi), \tilde{g}_k(\cdot, \alpha, \psi)\bigr) \leq
O_p \biggl(\frac
{S_{K}}{\delta
_{K}n^{1/2}} \biggr) = o_p(1).
\]
\upqed\end{pf*}


\begin{pf*}{Proof of Theorem~\ref{thmm: choice_K}}
We begin by placing the following restrictions on the sequence $p_n$:
\begin{longlist}[(ii)]
\item[(i)]$p_n \uparrow1$ and
\item[(ii)] for large $n$, $p_n \neq V_KV_\infty ^{-1}$ for any $K$.
\end{longlist}
Furthermore, the corresponding sequence $K^\ast$ must satisfy the
assumption of Lemma~\ref{lma: fnt_trnc}. Set $\epsilon _K =
\epsilon _K(n)
= |V_KV_\infty ^{-1}- p_n|$, $K = 1, \ldots, K^\ast$, where $K^\ast
$ is
given in~(\ref{eq: choice_K}), and define $\pi_{K^\ast} = \min\{
\epsilon _1,
\ldots, \epsilon _{K^\ast}\}$. Letting $S_{K^\ast}$ be defined as in
Lemma~\ref{lma: fnt_trnc} and ${\beta_{K^\ast} = n^{-1/2}
(S_{K^\ast
}\sum_{k = 1}^{K^\ast} \delta_k^{-1} )}$, we also require that
%
\begin{equation}
\label{eq: pi_rate} \biggl( \biggl(\frac{K^\ast}{n} \biggr)^{1/2}+
\beta_{K^\ast} + \gamma _{n} \biggr)\pi_{K^\ast}^{-1}
\rightarrow 0.
\end{equation}
None of these restrictions are contradictory.

Next, let $f_{i,K} = f_i(\cdot, K, \psi)$ and define
\[
\hat{V}_\infty = \frac{1}{n}\sum_{i=1}^nd(f_i,
f_\oplus )^2 \quad\mbox {and}\quad \hat{V}_K=
\hat{V}_\infty - \frac{1}{n}\sum_{i=1}^nd(f_i,
f_{i,K})^2.
\]
Observe that $\hat{V}_\infty - V_\infty = O_p(n^{-1/2})$ by the law
of large
numbers. Also, by (T3), for any $R > 0$,
\begin{eqnarray*}
P \Bigl(\max_{1 \leq K \leq K^\ast}\bigl\llvert (\hat{V}_\infty - \hat
{V}_K) - (V_\infty - V_K )\bigr\rrvert > R \Bigr)
&\leq&\frac{K^\ast}{R^2n} \max_{1 \leq K \leq K^\ast} E\bigl(d(f_1,
f_{1,K})^4\bigr)
\\
&=& O \biggl(\frac{K^\ast}{R^2n} \biggr).
\end{eqnarray*}
Hence,
\[
\max_{1 \leq K \leq K^\ast} \biggl\llvert \frac{\hat{V}_K}{\hat
{V}_\infty } -
\frac
{V_K}{V_\infty
}\biggr\rrvert = \max_{1 \leq K \leq K^\ast} \biggl\llvert
\frac{\hat
{V}_\infty - \hat{V}_K
}{\hat{V}_\infty
} - \frac{V_\infty - V_K}{V_\infty }\biggr\rrvert = O_p \biggl(
\biggl(\frac
{K^\ast
}{n} \biggr)^{1/2} \biggr).
\]

Define $\tilde{f}_{i,K} = \tilde{f}_i(\cdot, K, \psi)$. Then
observe that
\begin{eqnarray*}
\bigl\llvert (\hat{V}_\infty - \hat{V}_K) - (
\tilde{V}_\infty - \tilde {V}_K)\bigr\rrvert &\leq&
\frac
{1}{n}\sum_{i=1}^n\bigl\llvert
d(f_i, f_{i,K})^2 - d(f_i, \tilde
{f}_{i,K})^2\bigr\rrvert
\\
&\leq&\frac
{1}{n}\sum_{i=1}^nd(f_{i,K},
\tilde{f}_{i,K}) \bigl(2d(f_i, f_{i,K}) +
d(f_{i,K}, \tilde{f} _{i,K})\bigr).
\end{eqnarray*}
By using (T3), Lemma~\ref{lma: fnt_trnc} and the assumptions on the
sequence $K^\ast$, we find that
\[
\max_{1 \leq K \leq K^\ast}\bigl\llvert (\hat{V}_\infty -
\hat{V}_K) - (\tilde{V}_\infty - \tilde{V}_K )
\bigr\rrvert = O_p(\beta_{K^\ast}).
\]
By using similar arguments, we find that $\hat{V}_\infty - \tilde
{V}_\infty = O_p(\gamma _{n})$,
which yields
%
\begin{equation}
\label{eq: choice_K_rate} \max_{1 \leq K \leq K^\ast}\biggl\llvert \frac{V_K}{V_\infty } -
\frac
{\tilde{V}_K
}{\tilde{V}_\infty
}\biggr\rrvert = O_p \biggl( \biggl(
\frac{K^\ast}{n} \biggr)^{1/2}+ \beta _{K^\ast} + \gamma
_{n} \biggr).
\end{equation}

To finish, observe that, since $p_n \neq V_KV_\infty ^{-1}$ for any
$K$ when
$n$ is large, for such~$n$
\[
\bigl\{K^\ast\neq\tilde{K}^\ast \bigr\} = \biggl\{\max
_{1 \leq K
\leq
K^\ast
}\biggl\llvert \frac{V_K}{V_\infty } - \frac{\tilde{V}_K}{\tilde
{V}_\infty }
\biggr\rrvert > \pi _{K^\ast
} \biggr\}.
\]
Then, by (\ref{eq: choice_K_rate}), for any $\varepsilon > 0$ there
is $R > 0$
such that
\[
P \biggl(\max_{1 \leq K \leq K^\ast}\biggl\llvert \frac{V_K}{V_\infty } -
\frac{\tilde{V}_K
}{\tilde{V}_\infty }\biggr\rrvert > R \biggl( \biggl(\frac{K^\ast
}{n}
\biggr)^{1/2}+ \beta _{K^\ast} + \gamma _{n} \biggr)
\biggr) < \varepsilon
\]
for all $n$. Then, by (\ref{eq: pi_rate}), for $n$ large enough we have
$
P(K^\ast\neq\tilde{K}^\ast) < \varepsilon $.
\end{pf*}
\end{appendix}

\section*{Acknowledgments}
We wish to thank the Associate Editor and three referees for helpful
remarks that led to an improved version of the paper.




\begin{supplement}[id=suppA]
\stitle{The Wasserstein metric, Wasserstein--Fr\'echet mean,
simulation results and additional proofs}
\slink[doi]{10.1214/15-AOS1363SUPP} 
\sdatatype{.pdf}
\sfilename{aos1363\_supp.pdf}
\sdescription{The supplementary material includes additional
discussion on the Wasserstein distance 
and the rate of convergence of the Wasserstein--Fr\'echet mean is
derived. Additional simulation results are presented for FVE values
using the Wasserstein metric, similar to the boxplots in Figure~\ref
{fig: sim_fve}, which correspond to FVE values using the $L^2$ metric.
All assumptions are listed in one place. Lastly, additional proofs of
auxiliary results are provided.}
\end{supplement}

%


\printaddresses
\end{document}